\documentclass{elsarticle}
\usepackage{amsfonts}
\usepackage{fancyhdr}
\usepackage{comment}
\usepackage{amsmath}
\usepackage{changepage}
\usepackage{amssymb}
\usepackage{graphicx}%

\usepackage{microtype}

\usepackage{booktabs}
\usepackage{multirow}
\usepackage{subfig}

\usepackage{tikz}

\usepackage{pdfsync}
\usepackage[]{hyperref}
\hypersetup{colorlinks=false}

\usepackage{mathrsfs}
\usepackage{tensor}
\usepackage{bm}
\renewcommand{\theta}{\vartheta}
\renewcommand{\rho}{\varrho}
\renewcommand{\epsilon}{\varepsilon}
\newcommand{\sss}{\boldsymbol{s}}
\newcommand{\uu}{\boldsymbol{u}}
\newcommand{\xx}{\boldsymbol{x}}
\newcommand{\PP}{\boldsymbol{P}}

\newcommand{\CC}{\mathcal{C}}
\newcommand{\KK}{\mathrm{K}}

\newcommand{\QQ}{\mathrm{Q}}
\newcommand{\VV}{\mathrm{V}}

\newcommand{\Chris}{\Gamma}
\DeclareMathOperator{\LapB}{\Delta_{\mathrm{B}}}
\newcommand*{\de}{\mathop{}\!\mathrm{d}}

\newcommand{\Om}{\mathit{\Omega}}

\usepackage{color}

\newcommand{\dualityOmega}[1]{\left\langle #1 \right \rangle}
\newcommand{\dualityGamma}[1]{\left \langle \! \left\langle #1 \right \rangle\!\right\rangle}

\graphicspath{{./figures/}}

\begin{document}

\begin{frontmatter}
\title{NURBS-SEM: \\
a hybrid spectral element method on NURBS maps \\
for the solution of elliptic PDEs on surfaces}

\author[sissa]{Giuseppe Pitton}
\ead{gpitton@sissa.it}
\author[sissa]{Luca Heltai\corref{cor1}}
\ead{luca.heltai@sissa.it}
\cortext[cor1]{Corresponding author}
\address[sissa]{SISSA, Via Bonomea 265, Trieste, Italy}
\date{}
 
\begin{abstract}
Non Uniform Rational B-spline (NURBS) patches are a standard way to describe complex geometries in Computer Aided Design tools, and have gained a lot of popularity in recent years also for the approximation of partial differential equations, via the Isogeometric Analysis (IGA) paradigm. However, spectral accuracy in IGA is limited to relatively small NURBS patch degrees (roughly $p \leq 8$), since local condition numbers grow very rapidly for higher degrees. On the other hand, traditional Spectral Element Methods (SEM) guarantee spectral accuracy but often require complex and expensive meshing techniques, like transfinite mapping, that result anyway in \emph{inexact} geometries. In this work we propose a hybrid  NURBS-SEM approximation method that achieves spectral accuracy and maintains exact geometry representation by combining the advantages of IGA and SEM.

As a prototypical problem on non trivial geometries, we consider the Laplace--Beltrami and Allen--Cahn equations on a surface.
On these problems, we present a comparison of several instances of NURBS-SEM with the standard Galerkin and Collocation Isogeometric Analysis (IGA).

\end{abstract}


\end{frontmatter}

\section{Introduction}

Isogeometric analysis (IGA) is a numerical approximation paradigm for Partial Differential Equations (PDE) based on the exploitation of the same Non Uniform Rational B-spline (NURBS) basis functions used in Computer Aided Design (CAD) tools, and it has become a standard practice since its introduction in~\cite{IGAFirst,IGABook}.

It can be argued that the main merit of IGA is the removal of the mesh generation process required by other standard approximation methods (FEM, FVM, \dots) and the resulting tight interaction between CAD and analysis tools. Such a tight connection is particularly significant when the problem is formulated directly on lower dimensional manifolds embedded in the physical space, i.e., on surfaces in three dimensional space or curves in two or three-dimensional spaces, since in this case the Boundary Representation (BREP) standard which is common to all major CAD tools can be used \emph{as is}, without the need to create volumetric NURBS representations. This is the case for plates and shell structures~\cite{kiendl_isogeometric_2009, kiendl_isogeometric_2015}, Reissner--Mindlin shells \cite{benson_isogeometric_2010,dornisch_isogeometric_2013,dornisch_treatment_2014,uhm_tspline_2009}, boundary integral formulations~\cite{HeltaiKiendlDeSimone-2016-a,Heltai2014}, or when the physical problem of interest is naturally defined on manifolds~\cite{Bartezzaghi2015446,Bartezzaghi2016625,Dede2015} embedded in a Euclidean space.

Standard Finite Element approximation of such PDEs typically resort to an approximation of the geometry by means of piecewise polynomial patches.
This results in an approximation of the curvature of surfaces that may significantly affect the approximation quality. On the other hand, most of the geometries of practical interest can be represented exactly by B-splines or NURBS~\cite{Piegl}, making IGA an ideal candidate for these types of problems. 

NURBS-based Isogeometric analysis is a higher-order approach that allows an alternative refinement strategy with respect to the standard $h$- and $p$-refinements used in finite element analysis: the so called $k$-refinement~\cite{IGAFirst}, where the regularity of the basis functions is raised at each refinement stage.
This is of particular importance for eigenvalue problems, where IGA provides good approximations of a much larger part of the spectrum than standard FEA, as shown by~\cite{Cottrell2006} in a structural analysis context.
The accuracy and usability of $k$-refinement, however, is limited by cross-patch regularity.
For multiple-patch geometries, the global condition number of the system matrices presents the same behaviour under $p$- and $k$-refinement. The analysis of these methods~\cite{Sangallihpk,BAZILEVS2006b} proves that spectral accuracy can be achieved, but the error constants depend exponentially on the degree of the NURBS basis functions.
In practical applications this obstructs IGA from reaching spectral accuracy, which proves to be very efficient for relatively low degrees, and deteriorates when the basis functions' degree exceeds 8 or 9.

In this paper we present an alternative approximation strategy inspired by the Spectral Element Method (SEM), that aims at combining the advantages of IGA and SEM, similarly to what is done in~\cite{NME:NME2311} between FEM and NURBS. This approach is particularly efficient for those cases where it is important to treat the geometry exactly, but it would be desirable to achieve spectral accuracy on the solution of the PDE.
This is possible by breaking the iso-parametric paradigm of standard SEM and IGA, and constructing the push-forward of SEM basis functions through the NURBS description of the geometry. 

This article is organized as follows. In Section~\ref{sec:laplace_beltrami}, we briefly introduce the problem of solving an elliptic, second order PDE on a surface, and the Galerkin and Collocation strategies are introduced.
In Section~\ref{sec:methods}, we give a detailed description of the Spectral Element schemes introduced in this work, as well as a unified description encompassing some known methods such as IGA. The algebraic details leading to an efficient implementation of the Spectral Element Methods are also introduced in Section~\ref{sec:methods}, followed by a discussion on the implementation of essential and mixed boundary conditions.
In Section~\ref{sec:numerical_results}, we assess the performance of several hybrid NURBS-SEM instances by making a comparison with standard Galerkin and collocation Isogeometric Analysis for the approximate solution of the Laplace--Beltrami and Allen--Cahn equations on a surface, and discuss the results.
The Conclusions in Section~\ref{sec:conclusions} close the paper.

%
%
%
%
%
%

\section{Laplace--Beltrami and related Equations}
\label{sec:laplace_beltrami}

The Laplace--Beltrami operator can be thought of as an extension of the Laplace operator from a domain in $\mathbb{R}^n$ to a manifold. We briefly outline its derivation on a general Riemannian Manifold, together with its weak formulation. The subject is classical, and we refer for instance to~\cite{Lee,Klingenberg} for an introduction to Riemannian Geometry, and to~\cite{Lablee},~\cite{Rosenberg} for a detailed treatment on the properties of the Laplace--Beltrami operator (notably, its spectral properties).

Consider a bounded, orientable surface $\Om$ embedded in $\mathbb{R}^3$. Any regular surface can be described as a Riemannian Manifold, with metric tensor $g$. We assume that for each point $p \in \Om$, there exists a local coordinate description of $\Om$ that consists of a pair $(U,\xx^{-1})$, where $U$ is a neighbourhood of $p$ homeomorphic to an open set of $\mathbb{R}^2_+$ (namely, the upper half plane) and $\xx:\mathbb{R}^2\to U$ is a differentiable homeomorphism. Since we are mainly concerned with CAD applications, we consider the special case where $\xx$ is a tensor product of NURBS functions, that will be described in detail in Section~\ref{sec:surface_rep}.

In this representation, the local coordinates are typically given as maps from the reference domain $\widehat{\Om}$, that we take as the unit square, to the Euclidean space $\mathbb{R}^3$, i.e., $\sss=s^1 \mathbf{e}_1 +s^2\mathbf{e}_2 \in \widehat{\Om}$, and the surface is described through the coordinates of its points $\xx(\sss)=x^i(\sss) \mathbf{e}_i \subset\mathbb{R}^3$ as $\sss$ varies in $\widehat{\Om}$.

A basis for the tangent space at a point is given by the partial derivatives of $\xx$ with respect to the reference coordinates:
\begin{equation}
  \mathbf{g}_\alpha(\sss) := \frac{\partial x^i}{\partial s^\alpha}(\sss) \mathbf{e}_i,
  \label{eq:x_derivatives}
\end{equation}
where the summation with respect to $i$ is implied. Here and in the following, we adopt Einstein's summation convention on repeated indices. The metric tensor has the following local representation:
\begin{equation}
  g_{\alpha\beta}:= \mathbf{g}_\alpha \cdot \mathbf{g}_\beta = \xx_{,\alpha}\cdot\xx_{,\beta},
  \label{eq:def_g_ab}
\end{equation}
while the components of the inverse metric tensor $g^{\alpha\beta}$ are obtained by inversion of the $2\times2$ matrix $g_{\alpha\beta}$. We adopt the Einstein summation convention, with greek indices to indicate coordinates in $\widehat{\Om} \subset\mathbb{R}^2$, and latin indices for coordinates in $\mathbb{R}^3$.

Any Riemannian Manifold is naturally endowed with a Levi--Civita connection $\nabla$, that introduces a precise definition for the covariant derivatives of scalar functions or vector fields on a manifold.
The surface gradient of a scalar function $f$ at a point on $\Om$ is given by a differential form, whose representation in the local coordinate system is:
\begin{equation}
  \label{eq:surf_grad}
  \nabla f = \nabla_\mu f \mathbf{g}^\mu := \partial_\mu f \mathbf{g}^\mu  = \frac{\partial f}{\partial s^{\mu}}  \mathbf{g}^\mu,
\end{equation}
where $\{\mathbf{g}^\nu = g^{\nu \mu}\mathbf{g}_\mu\}$ is the dual basis to $\{\mathbf{g}_\mu\}$, i.e., $\mathbf{g}^\mu\cdot \mathbf{g}_\nu = \tensor{\delta}{^\mu_\nu}$ with $\delta$, the Kronecker delta, equal to one if $\nu = \mu$, and zero otherwise.
A representation of $\nabla f$ in the tangent space $T\Om$ is the vector field obtained by raising indices:
\begin{equation}
  \nabla f = \nabla^{\mu} f\mathbf{g}_\mu = g^{\mu\nu} \frac{\partial f}{\partial s^\nu}\mathbf{g}_\mu.
  \label{eq:nabla_up}
\end{equation}

For a vector field $\uu\in T\Om$, the covariant derivative is defined as:
\begin{equation}
  \nabla \uu:=\nabla (u^\alpha \mathbf{g}_\alpha) = \left(\frac{\partial u^\alpha}{\partial s^\beta}+\tensor{\Chris}{^\alpha_\beta_\gamma}u^\gamma\right)\mathbf{g}_\alpha\otimes\mathbf{g}^\beta,
  \label{eq:def_grad}
\end{equation}
where $\tensor{\Chris}{^\mu_\alpha_\beta}$ are the Christoffel symbols of the second kind, that can be expressed in terms of  partial derivatives of the metric tensor, as follows. First, the Christoffel symbols of the first kind are defined:
\begin{equation}
  \tensor{\Chris}{_\lambda_\mu_\nu}=\frac12\left(\partial_\nu g_{\lambda\mu}+\partial_\mu g_{\lambda\nu}-\partial_\lambda g_{\mu\nu}\right),
  \label{eq:coord_Chris1}
\end{equation}
then, the first index is raised:
\begin{equation}
  \tensor{\Chris}{^\lambda_\mu_\nu}=\tensor{g}{^\lambda^\rho}\tensor{\Chris}{_\rho_\mu_\nu}.
  \label{eq:coord_Chris2}
\end{equation}
With these notions, it is possible to introduce the Laplace--Beltrami operator acting on a scalar function as:
\begin{equation}
  \begin{aligned}
  \LapB :&= \nabla^\alpha\nabla_\alpha=g^{\alpha\beta}\nabla_\beta\nabla_\alpha 
         = g^{\alpha\beta}\nabla_\beta(\partial_\alpha) \\
        &= g^{\alpha\beta}\left[\partial^2_{\alpha\beta}+\tensor{\Chris}{^\mu_\alpha_\beta}\partial_\mu\right],
  \end{aligned}
  \label{eq:def_LapB}
\end{equation}

The derivatives of a function $u(\xx)$ defined on the surface can be obtained explicitly by the chain rule:
\begin{equation}
  \partial_\alpha u(\xx(\sss))=\frac{\partial u}{\partial x^i}\Big|_{\xx(\sss)}\frac{\partial x^i}{\partial s^\alpha}\Big|_{\sss}.
  \label{eq:def_u_der}
\end{equation}
From this point on, with some abuse of notation we will identify the surface $\Om$ with its parametrization $\xx(\sss)$.

\subsection{Laplace--Beltrami Equation}
\label{sec:laplace_beltrami_eq}

Let us now define the Laplace--Beltrami equation on a surface with boundary $\Om$.
Let $\CC^1(\Om)$ be the set of differentiable functions on $\Om$ with values in $\mathbb{R}$, and $\CC^2(\Om)$ the set of functions on $\Om$ with continuous second derivative, then the Laplace--Beltrami problem is: given a function $f\in \CC^0(\Om)$, find a function $u\in \CC^2(\Om)$ such that:
\begin{equation}
  \begin{cases}
  -\LapB u = f\qquad&\text{in }\Om \\
  u = h_D      &\text{on }\partial_D \Om \\
  \partial_nu = h_N      &\text{on }\partial_N \Om, \\
\end{cases}
  \label{eq:LaplaceBeltramiEq}
\end{equation}
where $\partial_D \Om$ and $\partial_N \Om$ form a partition of the boundary $\partial \Om$ of $\Om$, and $h_D:\partial_D \Om\to\mathbb{R}$ and $h_N:\partial_N\Om\to\mathbb{R}$ are respectively the given Dirichlet and Neumann boundary data.
For surfaces without boundary, the treatment is conceptually similar, except that boundary conditions cannot be imposed, and instead the average of $u$, or its value at a point must be specified.


In the following, we refer mainly to the variational form of Equation~\eqref{eq:LaplaceBeltramiEq}. Let $f\in H^{-1}(\Om)$, $h_D\in H^{\frac12}(\partial_D \Om)$ and $h_N\in H^{-\frac12}(\partial_N\Om)$ be given. We seek for a solution $u\in H^1_{D}(\Om)$ such that:
\begin{equation}
 (\nabla v,\nabla u)=  \dualityGamma{ v,h_N} + \dualityOmega{v,f}\qquad\forall v\in H^1_{0}(\Om),
  \label{eq:variational_LB}
\end{equation}
where the round brackets denote the $L^2(\Om)$ inner product, the brackets $\dualityGamma{\cdot,\cdot}$ are used to express the duality pairing between the Sobolev spaces $H^{-\frac12}(\partial_N\Om)$ and $H^{\frac12}(\partial_N\Om)$, while the brackets $\dualityOmega{\cdot,\cdot}$ are used for the duality pairing between the Sobolev spaces $H^{-1}(\Omega)$, and $H^1_D(\Omega)$ where
\begin{equation}
  H^1_{D}:=\{v\in H^1(\Om):\gamma_D u=h_D\},
  \label{eq:def_H1}
\end{equation}
for a suitable trace operator $\gamma_D:H^1_{D}(\Om)\to H^{\frac12}(\partial_D\Om)$.

Approximation methods for the variational problem~\eqref{eq:variational_LB} based on the Petrov--Galerkin projection stem from the introduction of two finite dimensional spaces $V^n,W^n$. An approximation $u^n\in V^n$ is defined as the solution to:
\begin{equation}
 (\nabla v,\nabla u^n) =   \dualityGamma{ v,h_N} + \dualityOmega{v,f}\qquad\forall v\in W^n.
  \label{eq:variational_LB_N}
\end{equation}
Let $\{\phi_j\}_{j=1}^n$ and $\{\psi_i\}_{i=1}^n$ be two sets of basis functions respectively for $V^n$ and $W^n$.
Then, the approximate solution can be expressed as $u^n(\xx)=\sum_{j=1}^n u^n_j\phi_j(\xx)$, and replacing this expansion in Equation~\eqref{eq:variational_LB_N} and choosing as test functions $v=\psi_i$ from a set of basis functions for $W^n$, Equation~\eqref{eq:variational_LB_N} is equivalent to the following linear system:
\begin{equation}
  \sum_{j=1}^nu^n_j(\nabla\psi_i,\nabla\phi_j)=\dualityGamma{ \psi_i,h_N } + \dualityOmega{\psi_i,f}\qquad \forall \psi_i,\, i=1,\dots,n.
  \label{eq:variational_algebra}
\end{equation}

Petrov--Galerkin schemes allow as well to construct approximations for nonlinear Laplace--Beltrami analogues. We discuss the nonlinear case by making the example of the Allen--Cahn equation:
\begin{equation}
  -\LapB u -u + u^3 = f
  \label{eq:allen_cahn}
\end{equation}
with the same boundary conditions discussed after Equation~\eqref{eq:LaplaceBeltramiEq}.
Before applying the Petrov--Galerkin projection, it is convenient to set up an iterative solution scheme. Here, we consider the following fixed-point iteration:
\begin{enumerate}
  \item an initial guess $u_0\in V^n$ is found by solving the linear part of the equation:
\begin{equation}
  -\LapB u_0 -u_0 = f;
  \label{eq:cahn_0}
\end{equation}
\item at the $n$-th iteration, a new approximation $u_{n+1}\in V^n$ is obtained by solving the linear equation:
\begin{equation}
  -\LapB u_{n+1} -u_{n+1} + u_n^2 u_{n+1}=f,
  \label{eq:cahn_n_1}
\end{equation}
whose weak form is: find $u_{n+1}\in V^n$ s.t.:
\begin{equation}
  (\nabla v,\nabla u_{n+1}) - (v,u_{n+1}) + (v,u_n^2u_{n+1}) = \dualityGamma{ v,h_N}  + \dualityOmega{v,f} \qquad \forall v\in W^n.
  \label{eq:cahn_n_weak}
\end{equation}
\item to judge on the quality of the current approximation $u_{n+1}$, it is convenient to consider the increment:
\begin{equation}
  d_n = \|u_{n+1}-u_n\|;
  \label{eq:def_increment}
\end{equation}
the iterative scheme is halted whenever $d_n$ is sufficiently small, e.g. by requiring that its order of magnitude be comparable with the machine epsilon, the smallest number representable in floating point arithmetic. 
\end{enumerate}

\subsection{Algebraic problem}
The following variational forms are naturally associated with the Laplace--Beltrami problems described in Section~\ref{sec:laplace_beltrami_eq}:
\begin{equation}
  \mathcal{A}(v,u) = (\nabla v,\nabla u), \qquad
  \mathcal{F}(v) = \dualityGamma{ v,h_N } +\dualityOmega{v,f} , \qquad
  \mathcal{N}(v,w,u) = (v,w^2u) - (v,u).
  \label{eq:varfsLB}
\end{equation}
The Laplace--Beltrami Equation~\eqref{eq:variational_LB_N} can be written in terms of these differential forms as:
\begin{equation}
  \mathcal{A}(v,u) = \mathcal F(v) \qquad \forall\; v \in W^n,
  \label{eq:AuvF}
\end{equation}
and the $n$-th iteration of the fixed-point scheme for the Allen--Cahn problem as:
\begin{equation}
  \mathcal{A}(v,u_{n+1}) + \mathcal{N}(v,u_n,u_{n+1}) = \mathcal{F}(v) \qquad \forall\; v \in W^n.
  \label{eq:varfANF}
\end{equation}
Evaluating the differential forms defined in Equation~\eqref{eq:varfsLB} on the basis sets $\{\phi_j\},\{\psi_i\}$ yields the matrices:
\begin{equation}
  \KK_{ij} = \mathcal{A}(\psi_i,\phi_j)
  \label{eq:KK_LB}
\end{equation}
for Laplace--Beltrami, and
\begin{equation}
  \KK_{ij} = \mathcal{A}(\psi_i,\phi_j) + \mathcal{N}(\psi_i,u_n,\phi_j)
  \label{eq:KK_AC}
\end{equation}
for a step of the Allen--Cahn fixed-point iteration. Note that in this latter case the matrix $\KK$ is not symmetric.
Associated with $\mathcal{F}$, the following column vector arises:
\begin{equation}
  \mathbf{f}_i = \mathcal{F}(\psi_i),
  \label{eq:F_psi}
\end{equation}
and it is also useful to introduce the column vector $\mathbf{u}$, whose components are the solutions' expansion coefficients: $\mathbf{u}_j = u^n_j$.
As a result, the Petrov--Galerkin approximation of Equations~\eqref{eq:variational_LB_N} and~\eqref{eq:cahn_n_1} is the solution of the following linear algebra problem:
\begin{equation}
  \KK\mathbf{u}=\mathbf{f}.
  \label{eq:weak_Kuf}
\end{equation}

An important class of Petrov--Galerkin discretizations is obtained from Equation~\eqref{eq:AuvF} with the choice of test functions $v=\psi_i=\delta(\xx_i)$, Dirac deltas centered at $n$ points $\{\xx_i\}_{i=1}^n\subset\Om$, where the space $V^n$ is at least twice differentiable.

In this case, the weak form expressed in equation~\eqref{eq:variational_LB} contains the surface gradient of Dirac deltas, which can only be interpreted in the usual weak sense, by integrating back by parts, and evaluating the result (corresponding to the original strong form of the problem) at the point $\xx_i$.

The resulting projection method is called a \emph{Collocation Method}, that strongly enforces the original differential equation at the collocation points $\{\xx_i\}_{i=1}^n$.
The linear systems arising from Collocation Methods have the same form of the linear system of Equation~\eqref{eq:weak_Kuf}, but with the following stiffness matrices:
\begin{equation}
  \KK_{ij}=-\LapB\phi_j(\xx_i), \qquad
  \KK_{ij} = -\LapB\phi_j(\xx_i) - \phi_j(\xx_i) + u_n^2(\xx_i)\phi_j(\xx_i),
  \label{eq:weak_def_Kuf}
\end{equation}
respectively for the linear and nonlinear case. In both cases the right hand side is just the evaluation of the forcing term at the collocation points:
\begin{equation}
  \mathbf{f}_i = f(\xx_i).
  \label{eq:weak_f_i}
\end{equation}
If the trial and test space coincide, $V^n=W^n$, the resulting projection is called a \emph{Galerkin projection}. When necessary, we will use the subscript $c$ to indicate \emph{collocation} matrices and the subscript $g$ to indicate \emph{Galerkin} matrices.



\section{Spectral Methods}
\label{sec:methods}

In this section, we introduce several instances of the Spectral Element Method, whose common feature is that of describing the computational domain exactly, through the surface NURBS representation coming e.g. from a Computer-Aided-Design (CAD) software. Two members of the family considered in the following coincide with the standard Isogeometric Analysis method and with an instance of the Isogeometric Collocation method. However, we regard that considering them as members of a large family of Spectral Element Methods is profitable in terms of thought and presentation economy.

\subsection{Surface representation}
\label{sec:surface_rep}
The most popular surface representation in the CAD community is that of Non-Rational Uniform B-splines (NURBS). Referring to~\cite{Piegl} for an introduction and thorough treatment on this subject, we briefly recall some important definitions. Here we focus only on tensor-product NURBS surfaces, wich are a subclass of the much more general family of NURBS parametrizations. Subdivision surfaces~\cite{NME:NME872} or Powell--Sabin B-splines~\cite{Manni2007,Speleers2012} are two possibile strategies when non-tensor product surfaces are required.

A NURBS description of a tensor-product bi-variate surface consists of two nondecreasing knot vectors, $\Theta^i=\{k^1=0,k^2,\dots,k^{n_i+p_i},k^{n_i+p_i+1}=1\}$, two positive weight vectors $W^i=\{w^1,\dots,w^{n_i}\}$ for $i=1,2$, and a set of $n=n_1n_2$ points $\{\PP^{i,j}\}_{i,j=1}^{n_1,n_2}\subset\mathbb{R}^3$, with $n_i \ge p_i+1$. Here it is assumed that the knot vectors $\Theta^i$ are \emph{open knot} vectors, meaning that the first and last knots are repeated $p_i+1$ times:
\[
  k^1 = k^2 = \dots = k^{p_i+1} = 0 \qquad \text{and} \qquad
  k^{n_i+1} = k^{n_i+2} = \dots = k^{n_i+p_i+1} = 1.
\]

The first step in the NURBS parametrization consists in defining the $n_d$ B-splines of degree $p$ with respect to the knot vector $\Theta^d$, $\{B^{(i,p)}\}_{i=1}^{n_d}$. The B-splines of degree zero are defined by:
\begin{equation}
  B^{(i,0)}(s)=
  \begin{cases}
    1\qquad &\text{if } k^i\le s<k^{i+1}\\
    0       &\text{otherwise}
  \end{cases}
  \label{eq:def_B0}
\end{equation}
and then the B-splines of degree $p$ are constructed by recursion using a convex combination of the B-splines of the previous degree, as described by the Cox--de Boor formula:
\begin{equation}
  B^{(i,p)}(s):=\frac {s-k^{i}}{k^{i+p}-k^{i}}B^{(i,p-1)}(s)+\frac {k^{i+p+1}-s}{k^{i+p+1}-k^{i+1}}B^{(i+1,p-1)}(s).
  \label{eq:def_Bip}
\end{equation}

The one-dimensional NURBS basis functions are defined by:
\begin{equation}
  N^{(i,p)}(s)=\frac{w^i B^{(i,p)}(s)}{\sum_{j=1}^{n_d}w^jB^{(j,p)}(s)},
  \label{eq:def_Ni}
\end{equation}
with $p=p_1,p_2$ depending on $d=1,2$.
The same procedure is repeated for the construction of $n_2$ splines and NURBS basis functions of degree $p_2$ in the other direction.

Finally, the surface is parametrized by weighting the grid points' coordinates with a tensor product of NURBS functions:
\begin{equation}
  \xx(\sss):=\sum_{i=1}^{n_1}\sum_{j=1}^{n_2}\PP^{i,j} N^{(i,p_1)}(s^1) N^{(j,p_2)}(s^2).
  \label{eq:def_xx_NURBS}
\end{equation}

\subsection{Approximation spaces}
\label{sec:basis}
For practical application of the Petrov--Galerkin projection, it is necessary to explicitly construct a basis set for both the approximation spaces $V^n$ and $W^n$. Here we consider only functions on the surface that are tensor product of 1D functions $\widehat{\phi}_i(s)$ defined on the interval $[0,1]$, composed with the inverse surface parametrization $\xx^{-1}$.

Consequently, the basis functions have the following shape:
\begin{equation}
  \phi_{\bm{i}}(\xx(\sss)):=\widehat{\phi}_{i_1}(s^1)\widehat{\phi}_{i_2}(s^2) \qquad \sss \in \widehat{\Om},
  \label{eq:def_phi_i}
\end{equation}
where $\boldsymbol{i}=(i_1,i_2)$ is an element of the index set:
\begin{equation}
  \mathcal{I}:=\{\boldsymbol{j}=(j_1,j_2), 1 \leq  j_1\le n_1, 1  \leq j_2 \leq n_2\}.
  \label{eq:def_I}
\end{equation}

In this work, we consider basis functions on the reference interval $[0,1]$, that originate different flavours of Petrov--Galerkin methods. The common features of the methods presented here are the following:
\begin{itemize}
  \item in all cases, the map $\xx$ between the reference element and the surface consists of a tensor product of NURBS functions;
  \item the basis functions are suitable in the sense of~\cite{Orszag} for a high-order method, meaning that the interpolation error for a smooth function converges to zero at an exponential rate as the degree of the basis functions is increased.
\end{itemize}

Two interesting choices for $\widehat{\phi}_i$ are B-spline and NURBS functions as defined in Equations~\eqref{eq:def_Bip} and~\eqref{eq:def_Ni} respectively.
Other two possibilities consist in Lagrange interpolants at Gauss--Lobatto points. Here we consider Gauss--Lobatto--Legendre (GLL) points, which are the $n-2$ zeros of the derivative of the Legendre Polynomial of degree $n-1$, plus the interval endpoints:
\begin{equation}
  s_j^{\mathrm{GLL}} = 
  \begin{cases}
    j\text{-th zero of }P_{n-1}'(s)\qquad &\text{for }j=1,\dots,n-2 \\
    0 & \text{for }j=n-1 \\
    1 & \text{for }j=n \\
  \end{cases}.
  \label{eq:def_s_gll}
\end{equation}
and Gauss--Lobatto--Chebyshev (GLC) points:
\begin{equation}
  s^{\mathrm{GLC}}_j = \cos \left(\pi + \frac{j-1}{n-1}\pi\right) \qquad \text{for } j=1,\dots,n.
  \label{eq:def_s_glc}
\end{equation}
Then, the Lagrange interpolants are defined as:
\begin{equation}
  L_j(s) = \prod_{\substack{ i=1 \\ i\neq j}}^n \frac{s-s_i^\mathrm{GLL}}{s_j^\mathrm{GLL}-s_i^\mathrm{GLL}} \qquad \text{for }j=1,\dots, n,
  \label{eq:def_lagrange_interp}
\end{equation}
for the Lagrange case, and
\begin{equation}
  T_j(s) = \prod_{\substack{ i=1 \\ i\neq j}}^n \frac{s-s_i^\mathrm{GLC}}{s_j^\mathrm{GLC}-s_i^\mathrm{GLC}} \qquad \text{for }j=1,\dots, n,
  \label{eq:def_chebyshev_interp}
\end{equation}
for the Chebyshev case.
In the following, we denote with $L^{n}_j$ the $j$-th Lagrange interpolant on $n$ GLL points, and with $T^n_j$ the $j$-th Lagrange interpolant on $n$ GLC points.

By taking different combinations of basis functions for the trial $V^n$ and test $W^n$ spaces, we generate the family of methods listed in Table~\ref{tab:hom_list}. The names of the different methods are chosen to be self-explanatory and compatible with the literature. For the cases where the trial and test space coincide, the method is classified as a Galerkin method, and any time the test space is formally spanned by Dirac deltas, the method is classified as a Collocation method.
\begin{table}
  \centering
  \begin{tabular}{llcc}
    \toprule
      Method name & Acronym & $V^n$ & $W^n$ \\
    \midrule
    B-spline Galerkin & SG & $B^{(i,p)}$ & $B^{(j,p)}$  \\
    B-spline Collocation & SC & $B^{(i,p)}$ & $\delta(\xx(s^\mathrm{Greville}_j))$   \\
    Isogeometric Galerkin & IG & $N^i$ & $N^j$ \\
    Isogeometric Collocation & IC & $N^i$& $\delta(\xx(s^\mathrm{Greville}_j))$  \\
    Chebyshev Collocation & CC& $T^n_i$ & $\delta(\xx(s^{\mathrm{GLC}}_j))$  \\
    Chebyshev Galerkin & CG & $T^n_i$ & $T^n_j$ \\
    Legendre Galerkin & LG & $L^n_i$ & $L^n_j$ \\
    \bottomrule
  \end{tabular}
  \caption{Catalogue of the 7 candidates for high-order methods considered in this work.}
  \label{tab:hom_list}
\end{table}


Below we briefly describe each numerical method that will be considered in the following. Our description is by no means exhaustive, and for a treatment of the first four methods of Table~\ref{tab:hom_list} we refer to~\cite{IGABook},~\cite{Schillinger2013} and~\cite{Auricchio2012}. However, we let the reader beware that none of these references is concerned with the application of Isogeometric Methods in the context of high-order refinement.
For a treatment of high-order methods similar, but not equivalent, to the last three rows of Table~\ref{tab:hom_list}, we refer to~\cite{Orszag},~\cite{CHQZ2},~\cite{CHQZ3}.

\emph{B-spline Galerkin Method (SG)}, where the B-spline basis functions defined in~\eqref{eq:def_Bip} are used as basis both for the trial and test spaces.

\emph{B-spline Collocation Method (SC)}. In this case, the B-spline basis functions of equation~\eqref{eq:def_Bip} are used as basis for the trial space, and Dirac deltas centered on the Greville abscissae of the B-spline functions are used as basis for the test space. 

\emph{Isogeometric Galerkin Method (IG)}. In this method, the NURBS basis functions defined in~\eqref{eq:def_Ni} are used as basis both for the trial and test spaces. A more consistent name for this method would be \emph{NURBS Galerkin Method}, since however it has been named \emph{Isogeometric Analysis} by its ideators in~\cite{IGAFirst}, we stick to this by now standard naming.

\emph{Isogeometric Collocation Method (IC)}. In this case, the NURBS basis functions~\eqref{eq:def_Ni} are used for constructing the trial space, and Dirac deltas centered on the Greville abscissae of the NURBS basis functions are used as basis for the test space.

\emph{Chebyshev Collocation Method (CC)}. For this method, the reference basis functions for the trial space are Lagrange interpolants on Gauss--Lobatto--Chebyshev nodes, and the test functions are Dirac deltas centered in the same Gauss--Lobatto--Chebyshev nodes.
The basis functions are constrained to be continuous together with their normal derivative across elements, as done in Chebyshev multipatch methods, see~\cite[p.~339]{CHQZ3} for details.

\emph{Chebyshev Spectral Element Method (CG)}. In this method, the basis functions for both trial and test spaces are Lagrange interpolants on Gauss--Lobatto--Chebyshev nodes. The Galerkin method in this case is based on the weak form~\eqref{eq:variational_LB}, with inner products weighted by the function
\begin{equation}
  w(x) = \frac{1}{\sqrt{1-x^2}}.
  \label{eq:def_weight_cheb}
\end{equation}
As in the previous case, inter-element continuity of the basis functions and their normal derivative is explicitly enforced.
One important characteristic of this and the following methods is that the quadrature formulas used for the evaluation of the integrals are based on the same nodes used to define the basis functions. This results in a ``variational crime'' due to under-integration, but also in a diagonal mass matrix, and allows for efficient integral evaluation. In some books, the methods with integration by quadratures are called ``SEM with numerical integration'', or ``SEM-NI''.

\emph{Legendre Spectral Element Method (LG)}. In this method, the basis functions for both trial and test spaces are Lagrange interpolants on Gauss--Lobatto--Legendre nodes. 
The basis functions are only continuous across contiguous elements.

\begin{figure}[tbp]
  \centering
\subfloat[]{
  \includegraphics[width=0.45\textwidth]{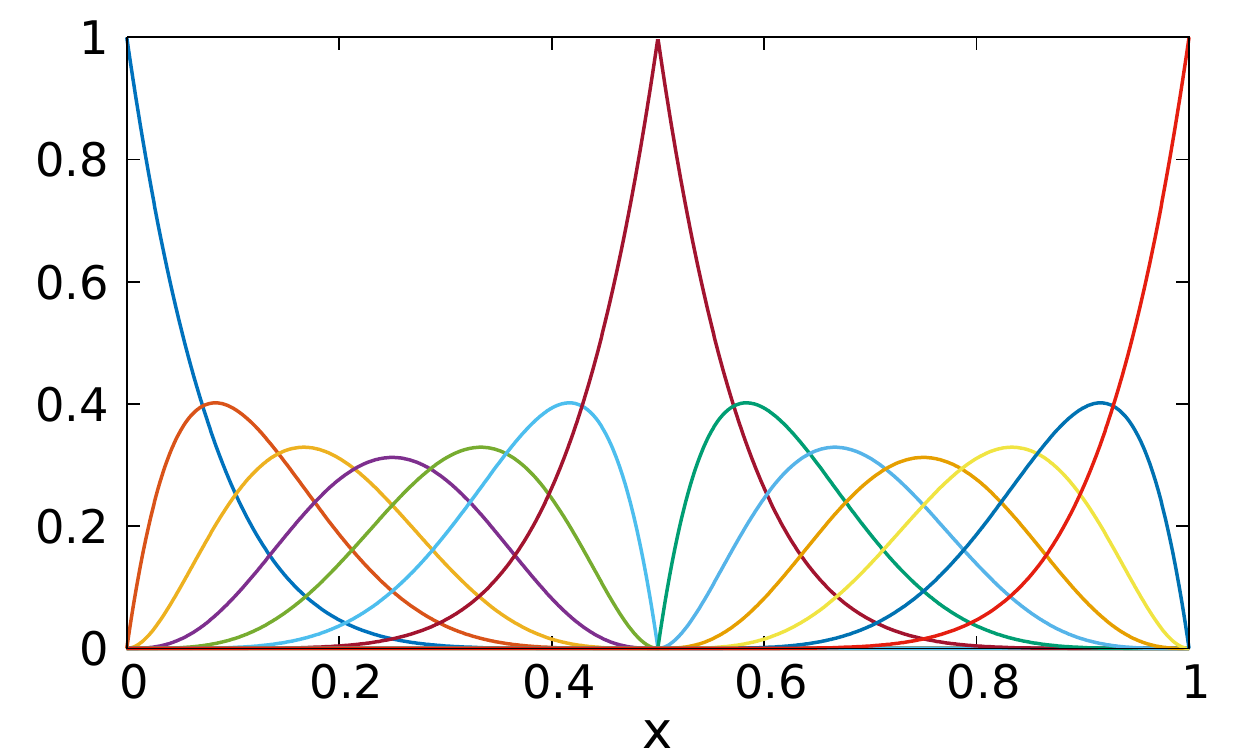}
}
\subfloat[]{
  \includegraphics[width=0.45\textwidth]{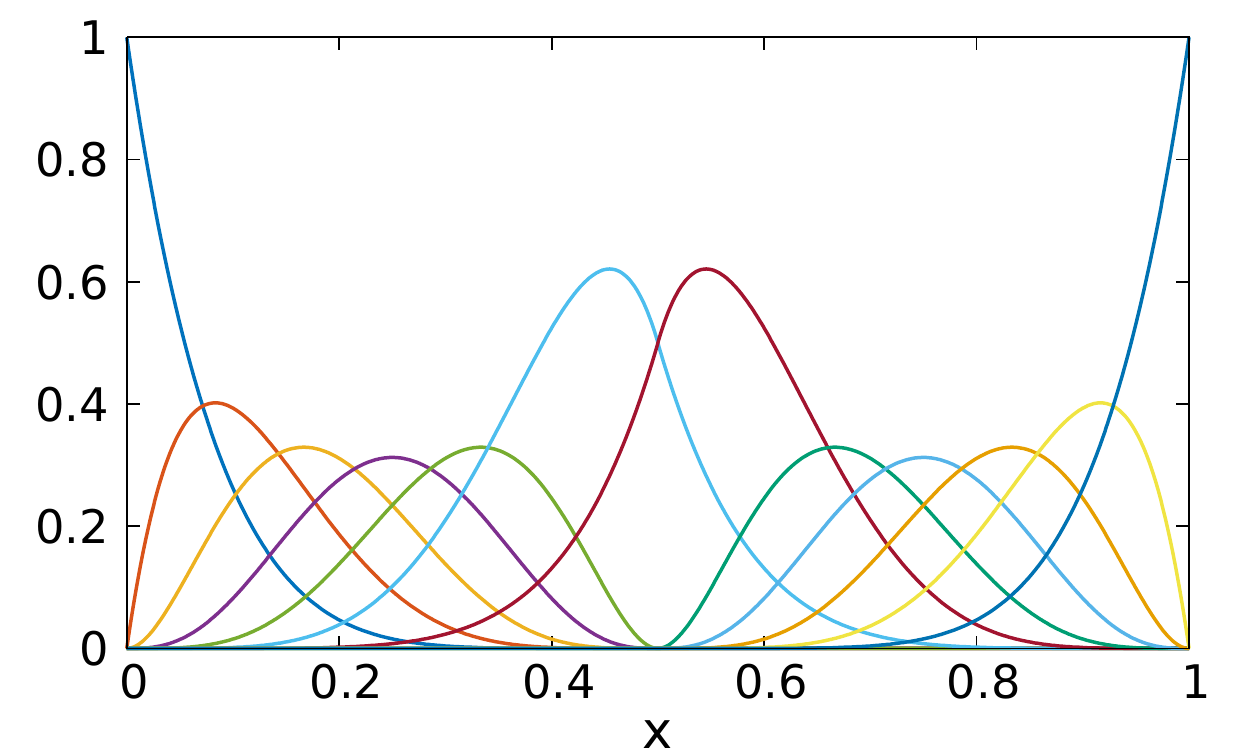}
} \\
\subfloat[]{
  \includegraphics[width=0.45\textwidth]{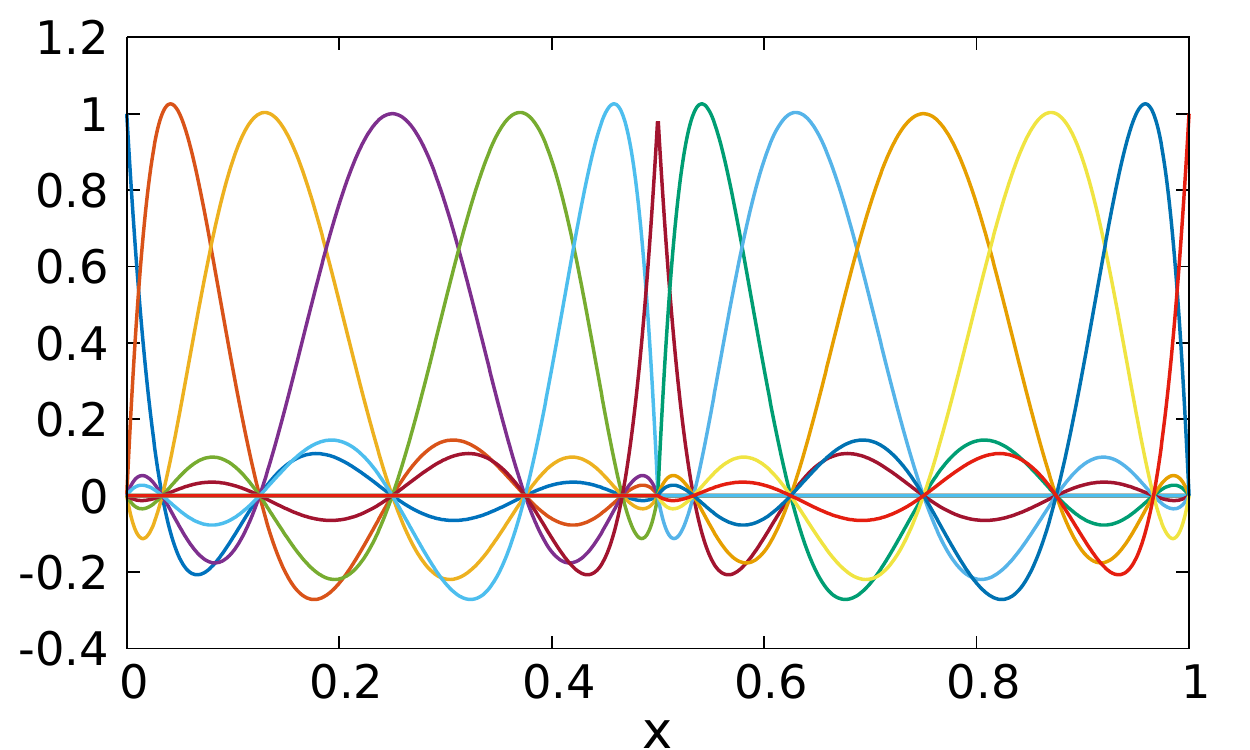}
}
\subfloat[]{
  \includegraphics[width=0.45\textwidth]{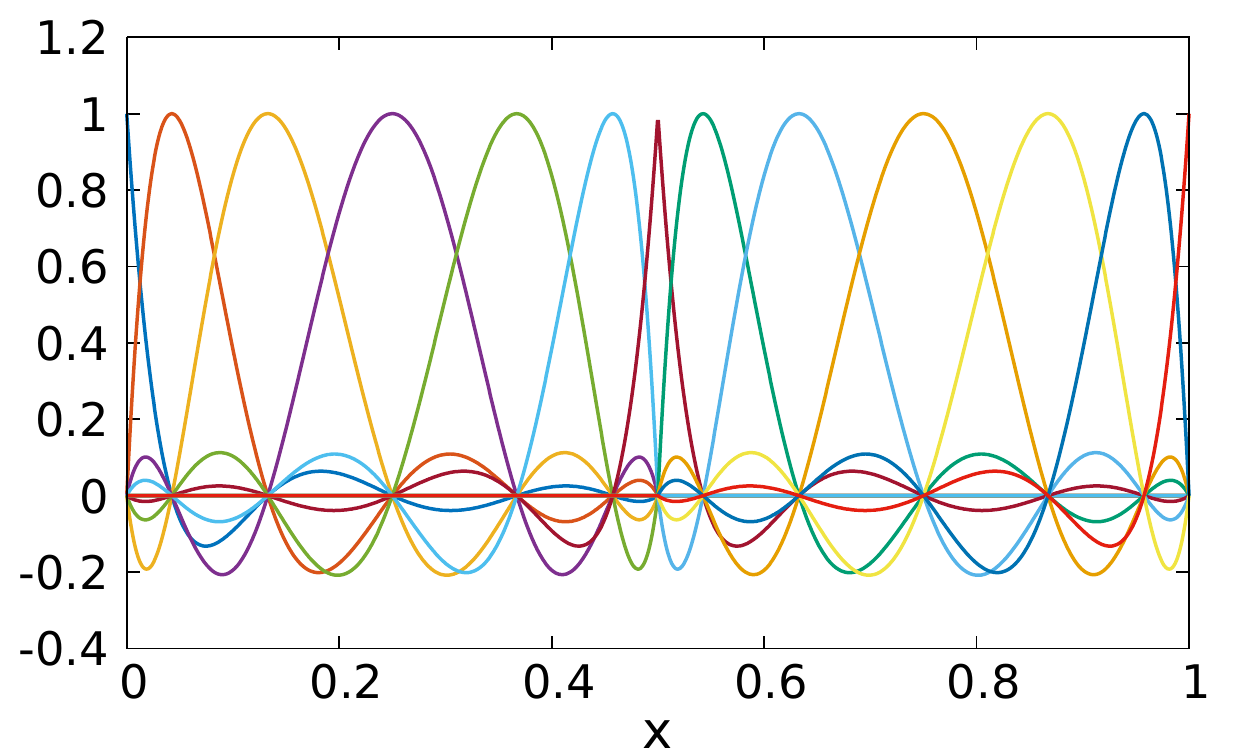}
}
  \caption{Basis functions of degree 7 on the reference interval for the case of two equal elements for $\CC^0$ NURBS (a), $\CC^1$ NURBS (b), Chebyshev Spectral Element (c), Legendre Spectral Element (d) methods.}
  \label{fig:basis_funs_1d}
\end{figure}

To provide some intuition over the different basis functions involved, we draw in Figure~\ref{fig:basis_funs_1d} the one-dimensional basis functions of degree 7 on the reference interval $[0,1]$ for the cases with 2 elements.
The B-spline functions are visually similar to the NURBS in many contexts.
In Figure~\ref{fig:basis_funs_1d} the two elements are the intervals $[0,0.5]$ and $[0.5,1]$, and in Figure~\ref{fig:basis_funs_1d} (b) the $\CC^1$ inter-element continuity of the NURBS functions is evident.

\subsection{Discrete formulation}
\label{sec:discrete_formulation}
Although the matrices that appear in the linear Equation~\eqref{eq:weak_Kuf} could in principle be computed naively as by their definition in Equations~\eqref{eq:KK_LB} and~\eqref{eq:KK_AC}, this often is not the most efficient way, especially for high-order methods.
In this section we give some details on how the discrete formulations of Section~\ref{sec:basis} could be set up and solved in a reasonably efficient manner.

Significant efficiency savings can be achieved by exploiting the tensor product form of the surface parametrization (Equation~\eqref{eq:def_xx_NURBS}) and of the basis function definition (Equation~\eqref{eq:def_phi_i}). 
These tensor product structures allow in turn to express the discrete operators in tensor-product form, that as we shall see drastically reduces the memory and the number of operations required to compute a matrix-vector product. 
This fact holds true both for methods based on the B-spline familty and for methods based on Lagrange interpolants.
For a detailed reference on the efficient implementation of tensor-product high-order methods, we refer to~\cite{DevilleFischerMund}.

The first step in the discretization of the variational form~\eqref{eq:variational_LB_N} is the evaluation of the inner products $(\nabla u,\nabla v)$ on all the basis functions. This operation can be performed by pulling back the integrals to the reference domain $\widehat{\Om}:=[0,1]^2$, as follows:
\begin{equation}
  \begin{aligned}
\int_{\Om} \nabla\phi_{\bm{i}}(\xx)\cdot\nabla\phi_{\bm{j}}(\xx)\de\xx
&= \sum_{k=1}^3 \int_{\Om}\frac{\partial \phi_{\bm{i}}(\xx)}{\partial x^k}\frac{\partial \phi_{\bm{j}}(\xx)}{\partial x^k} \de \xx \\
&= \sum_{k=1}^3\sum_{\alpha,\beta=1}^2\int_{\widehat{\Om}} \frac{\partial\widehat{\phi}_{\bm{i}}(\xx(\sss))}{\partial x^k}\frac{\partial x^k}{\partial s_\alpha} \frac{\partial\widehat{\phi}_{\bm{j}}(\xx(\sss))}{\partial x^k}\frac{\partial x^k}{\partial s_\beta} J(\sss)\de\sss 
=: K_{\bm{i}\bm{j}},
\end{aligned}
  \label{eq:def_int_x}
\end{equation}
where $J=\sqrt{|\det g_{\alpha\beta}|}$, and the the hat is used to denote quantities defined on the reference domain.
The Neumann boundary terms are also evaluated on the boundary of the reference domain, $\partial\widehat{\Om}_N$:
\begin{equation}
  \begin{aligned}
    \int_{\partial\Om_N} \phi_{\bm{i}}(\xx)\nabla\phi_{\bm{j}}(\xx) \cdot \bm{\nu} \de\xx
    &=  \int_{\partial\Om_N} \phi_{\bm{i}}(\xx) \sum_{k=1}^3\frac{\partial\phi_{\bm{j}}(\xx)}{\partial x^k} \nu^k \de\xx \\
    &=  \int_{\partial\widehat{\Om}_N} \widehat{\phi}_{\bm{i}}(\sss) h_N(\sss)J(\sss)\de\sss, \\
  \end{aligned}
  \label{eq:def_neumann_int}
\end{equation}
where $\widehat{\bm{\nu}}=(\widehat{\nu}^1,\widehat{\nu}^2)$ is the outer normal to the reference domain $\widehat{\Om}$, and $\widehat{\nu}$ is given in components by:
\begin{equation}
  \nu^k = \widehat{\nu}^\alpha\frac{\partial x^k}{\partial s^\alpha}.
  \label{eq:normal_components}
\end{equation}

The integrals on the right hand side of Equations~\eqref{eq:def_int_x} and~\eqref{eq:def_neumann_int} are then evaluated by quadratures. To this end, let $\{(\xi_k,\omega_k)\}_{k=1}^{Q_1}$ be the quadrature point-quadrature weight couples for the first coordinate, and $\{(\eta_l,\rho_l)\}_{l=1}^{Q_2}$ the quadrature points and weights for the second coordinate. Here $Q_1$ and $Q_2$ denote the number of quadrature points chosen for the two directions. It is convenient to introduce the auxiliary matrices:
\begin{equation}
  (G_g^{\alpha\beta})_{kl}=g^{\alpha\beta}(\xi_k,\eta_l)J(\xi_k,\eta_l)\omega_k\rho_l
  \label{eq:def_Gab}
\end{equation}
and
\begin{equation}
  (M)_{ki} = \widehat{\phi}_i(\xi_k) \qquad (D_1)_{ki} = \widehat{\phi}_i'(\xi_k).
  \label{eq:def_MD}
\end{equation}
The right hand side of equation~\eqref{eq:def_int_x} can be approximated by:
\begin{equation}
  K_g =
  \begin{pmatrix}
    M\otimes D_1 \\
    D_1\otimes M
  \end{pmatrix}^T
  \begin{bmatrix}
    G_g^{11} & G_g^{12} \\
    G_g^{21} & G_g^{22}
  \end{bmatrix}
  \begin{pmatrix}
    M\otimes D_1 \\
    D_1\otimes M
  \end{pmatrix}
  \label{eq:def_tensor_mat}
\end{equation}
which is a sequence of matrix-matrix products written in block tensor form.
The special block-tensor product structure of the expression~\eqref{eq:def_tensor_mat} allows for a reduction of the computational cost with respect to a general basis lacking the tensor product structure. Indeed, suppose that in both directions there are $N$ basis functions, and that the same $Q$-point quadrature rules are applied. Then, the number of operations required for the evaluation of the integral~\eqref{eq:def_int_x} decreases from $O(Q^2N^2)$ to $O(Q^2+2QN)$.

The Neumann boundary term in Equation~\eqref{eq:def_neumann_int} is computed similarly, and it represents a known term in the resulting algebraic system.

The action of the Laplace--Beltrami operator in the collocative case is even simpler, since it is sufficient to evaluate the derivatives of the basis functions and of the metric tensor on the collocation points, and then assemble point by point the coordinate expression of Equation~\eqref{eq:def_LapB}.
More precisely, let us introduce the collocation matrix for the second order derivatives of the basis functions:
\begin{equation}
  (D_2)_{ki} = \widehat{\phi}_i''(\xi_k),
  \label{eq:def_D2}
\end{equation}
we define the matrices $G_c^{\alpha\beta}$ for the collocation case as:
\begin{equation}
  (G_c^{\alpha\beta})_{kl}=g^{\alpha\beta}(\xi_k,\eta_l),
  \label{eq:def_Gab_collocation}
\end{equation}
and we introduce the pointwise evaluation of the Christoffel symbols:
\begin{equation}
  (\tensor{C}{^\mu_\alpha_\beta})_{kl}=g^{\alpha\beta}(\xi_k,\eta_l)\tensor{\Chris}{^\mu_\alpha_\beta}(\xi_k,\eta_l).
  \label{eq:def_C_collocation}
\end{equation}
The resulting discretization of the Laplace--Beltrami operator is:
\begin{multline}
  K_c =
  \begin{bmatrix}
    G_c^{11} & G_c^{12} \\
    G_c^{21} & G_c^{22}
  \end{bmatrix}
  \begin{bmatrix}
    M\otimes D_2 & D_1\otimes D_1 \\
    D_1\otimes D_1 & D_2\otimes M
  \end{bmatrix}
  +
  \begin{bmatrix}
    \tensor{C}{^1_1_1} & \tensor{C}{^2_1_1} \\
    \tensor{C}{^1_2_1} & \tensor{C}{^2_2_1} \\
  \end{bmatrix}
  \begin{pmatrix}
    M\otimes D_1 \\
    D_1\otimes M
  \end{pmatrix}
  + \\
  \begin{bmatrix}
    \tensor{C}{^1_1_2} & \tensor{C}{^2_1_2} \\
    \tensor{C}{^1_2_2} & \tensor{C}{^2_2_2} \\
  \end{bmatrix}
  \begin{pmatrix}
    M\otimes D_1 \\
    D_1\otimes M
  \end{pmatrix}.
  \label{eq:LapB_collocation}
\end{multline}
We remark that in Equation~\eqref{eq:LapB_collocation}, the collocation matrices $G_c^{\alpha\beta}$, defined according to Equation~\eqref{eq:def_Gab_collocation}, are different from the Galerkin matrices $G_g^{\alpha\beta}$, defined according to Equation~\eqref{eq:def_Gab}, 

For the Chebyshev-SEM method, the stiffness matrix is computed exactly as in Equation~\eqref{eq:LapB_collocation}, replacing $G_c^{\alpha\beta}$ with $G_g^{\alpha\beta}$ and replacing $\tensor{C}{^\mu_\alpha_\beta}$ with:
\begin{equation}
  (\tensor{S}{^\mu_\alpha_\beta})_{kl} = g^{\alpha\beta}(\xi_k,\eta_l)\tensor{\Chris}{^\mu_\alpha_\beta}(\xi_k,\eta_l)J(\xi_k,\eta_l)\omega_k\rho_l.
  \label{eq:def_C_integral}
\end{equation}

For both Galerkin and Collocation methods based on Chebyshev points, Neumann boundary conditions can be imposed by collocation. This can be achieved by replacing the rows related to the collocation points lying on $\partial\Om_N$ with:
\begin{equation}
  \sum_l \phi_{\bm{i}}(\xi_N,\eta_l)\nabla\widehat{\phi}_{\bm{j}}(\xi_N,\eta_l)\cdot\widehat{\bm{\nu}}(\xi_N,\eta_l),
  \label{eq:def_neumann_collocation}
\end{equation}
where to simplify the presentation we made the hypothesis that the Neumann boundary is located at the points with $s^1=\xi_N$ for all $s^2$, and the corresponding rows on the right hand side should be replaced by the known value $h_N(\xi_N,\eta_l)$.

For B-spline and Isogeometric methods, the imposition of Neumann boundary conditions by collocation is still an active research area. Since there are no simple, established solution to this problem, we do not discuss it here, and refer instead to~\cite{AURICCHIO2010a, DELORENZIS201521}.

Computing the system's matrix in the Allen--Cahn case requires one more step, namely the efficient evaluation of the integral coming from the fixed point linearization
\begin{equation}
  \int_\Om \phi_{\boldsymbol{i}}(\xx) u_n^2(\xx)\phi_{\boldsymbol{j}}(\xx)\de\xx =
  \int_{\widehat{\Om}} \widehat{\phi}_{\boldsymbol{i}}(\sss) u_n^2(\xx(\sss)) \widehat{\phi}_{\boldsymbol{j}}(\sss) J(\sss) \de\sss.
  \label{eq:int_NN}
\end{equation}
This additional term is discretized as:
\begin{equation}
  (\psi_{\boldsymbol{i}(i,j)},u_n^2\phi_{\boldsymbol{j}(k,l)}) = (M_{ri}\otimes M_{sj}) N_{rs} (M_{rk}\otimes M_{sl}),
  \label{eq:allen_nonlinearity}
\end{equation}
where the index sets $\boldsymbol{i}(i,j)$ and $\boldsymbol{j}(k,l)$ are introduced to pass from the two-dimensional matrix notation on the left hand side to the multi-dimensional indexing of the right hand side. On the right hand side the summation with respect to $r$ and $s$ is implied.
The nonlinearity $N$ can either be evaluated at the quadrature points as:
\begin{equation}
  N_{rs} = (\phi_{\boldsymbol{k}}(\xi_r,\eta_s)u_{n,\boldsymbol{k}})^2 J(\xi_r,\eta_s)\omega_r\rho_s
  \label{eq:allen_N_kl_galerkin}
\end{equation}
in the case of Galerkin methods, or at the collocation points as:
\begin{equation}
  N_{rs} = (\phi_{\boldsymbol{k}}(\xi_r,\eta_s)u_{n,\boldsymbol{k}})^2,
  \label{eq:allen_N_kl_collocation}
\end{equation}
for collocation methods.
In Equations~\eqref{eq:allen_N_kl_galerkin} and~\eqref{eq:allen_N_kl_collocation}, the symbol $\phi_{\boldsymbol{k}}u_{n,\boldsymbol{k}}$ is a shortcut for the sum:
\begin{equation}
  \sum_{i=1}^{n_1}\sum_{j=1}^{n_2}\phi_i\phi_j u_{n,ij}.
  \label{eq:def_summation_N}
\end{equation}
Note that in Equation~\eqref{eq:def_summation_N}, the index $n$ refers to the $n$-th fixed point iteration, while $i$ and $j$ sum over all the basis functions.

In the case of Galerkin methods, aliasing errors, if present, may be reduced by evaluating $N$ via higher-order quadrature rules. Aliasing, however, has not been an issue in the present work, since the resolution was sufficient to represent the nonlinearity with a good precision. For details on how aliasing may affect a computation not sufficiently resolved, see~\cite{Boyd}.

In the case of multi-element or multi-patch discretizations, the steps outlined above are repeated elementwise, and summed in the global stiffness matrix with the appropriate numbering of the degrees of freedom, as in standard Finite Element codes.

\subsection{Essential boundary conditions}
After treating natural boundary conditions in Section~\ref{sec:discrete_formulation}, we now discuss essential, or Dirichlet, boundary conditions.
For the methods based on Lagrange interpolants, the essential boundary conditions are imposed by row elimination.
The rows related to the boundary degrees of freedom are replaced by the corresponding rows of the identity matrix:
\begin{equation}
  \KK_{ij} = \delta_{ij} \qquad\forall\; i\text{ such that } \xx_i\in\partial_D\Omega,
  \label{eq:K_boundary}
\end{equation}
and at the right hand side, the Dirichlet datum is imposed at the corresponding degree of freedom:
\begin{equation}
  \mathbf{f}_i=h_D(\xx_i)  \qquad\forall\; i\text{ such that } \xx_i\in\partial_D\Omega.
  \label{eq:f_boundary}
\end{equation}

For the methods based on B-spline or NURBS basis functions, such a direct approach is not feasible since the basis is not interpolatory. We resort to a least-squares enforcement of essential boundary conditions.
The least-squares problem requires a set of points on the boundary $\{\overline{\xx}_i\}_{i=1}^q$ and an index-set $\bm{i}$ where the indices of the basis functions different from zero at the boundary are stored. The basis functions indexed by $\bm{i}$ are evaluated at the boundary points, forming the matrix $\VV\in\mathbb{R}^{q\times n_b}$:
\begin{equation}
  \VV_{ij} = \phi_{\bm{i}(j)}(\overline{\xx}_i) \qquad \text{for }j=1,\dots,n_b,\; i=1,\dots,q
  \label{eq:def_Vij}
\end{equation}
where $n_b$ is the cardinality of $\bm{i}$.
The least-squares problem is well-posed if and only if $q\ge n_b$.
Similarly, we introduce the array $\mathbf{q}\in\mathbb{R}^{q}$, containing the evaluation of the Dirichlet datum on the interpolation points:
\begin{equation}
  \mathbf{q}_j = h_D(\overline{\xx}_j) \qquad\text{for }j=1,\dots,q,
  \label{eq:def_qj}
\end{equation}
and the boundary restriction matrix $\QQ\in\mathbb{R}^{n_b\times n}$. We also introduce the following spaces:
\begin{equation}
  V^n_\circ =\{v\in V^n: \gamma_{\partial\Om}v=0\},
  \label{eq:def_VN_circ}
\end{equation}
\begin{equation}
  V^n_\partial = V^n \setminus V^n_\circ = \text{span} \{ \phi_{\bm{i}(j)}\}_{j=1}^{n_b},
  \label{eq:def_VN_partial}
\end{equation}
and we denote by $\mathbf{y}_\circ$ any vector of $\mathbb{R}^n$ containing the expansion coefficients in $\{\phi_j\}$ of any function with vanishing trace on the Dirichlet part of the boundary:
\begin{equation}
  (\mathbf{y}_\circ)_{\bm{i}(j)}=0 \qquad \text{for } j=1,\dots,n_b,
  \label{eq:def_y_circ}
\end{equation}
and with $\mathbf{y}_\partial$ we denote any vector of $\mathbb{R}^n$ containing the expansion coefficients of any function whose trace is not identically zero on $\partial_D\Om$, and whose interior coefficients vanish:
\begin{equation}
  (\mathbf{y}_\partial)_j= 0 \qquad \text{for } j \notin \{\bm{i}(k)\}_{k=1}^{n_b}.
  \label{eq:def_y_partial}
\end{equation}

Notice that functions in $V_\partial^n$ are in general not vanishing in the interior of the domain pointwise. They decay as we approach the interior of the domain in a mesh dependent way, following the decay of the B-spline basis functions whose value on the boundary is non-zero.

The restriction operator is then defined as:
\begin{equation}
  \begin{cases}
    \QQ \mathbf{y}_\circ = \mathbf{0} \\ 
    \QQ\QQ^T \mathbf{y}_{\partial}=\mathbf{y}_\partial \\ 
  \end{cases}
  \label{def:restrictionQ}
\end{equation}
We remark that $\QQ^T:\mathbb{R}^{n_b}\to\mathbb{R}^n$ takes a vector of boundary ``values'' and extends it to zero on the interior, and $\QQ^T\QQ\in\mathbb{R}^{n\times n}$ returns the boundary lifting of a vector.

The least squares imposition of Dirichlet boundary conditions is achieved through the augmented linear system:
\begin{equation}
  \begin{bmatrix}
    \KK & \QQ^T\VV^T\VV \\
    \VV^T\VV\QQ & 0
  \end{bmatrix}
  \begin{pmatrix}
    \mathbf{u} \\ \mathbf{\lambda}
  \end{pmatrix}
  =
  \begin{pmatrix}
    \mathbf{f} \\ \VV^T\mathbf{q}
  \end{pmatrix},
  \label{eq:LS_Algebra}
\end{equation}
where $\mathbf{\lambda}\in\mathbb{R}^{n_b}$ is a vector of auxiliary Lagrange multipliers.
A side-effect of this augmentation is the growth of the matrix condition number, but in our numerical experiments this is contained to a factor of $\simeq100$. The number of unknowns usually does not increase by more than $\simeq20\%$, but this depends on the boundary to surface ratio of the domain.
Alternatively, one could impose the boundary conditions by Nitsche method, as done, for example, in~\cite{RotundoKimJiang-2016-a}.

\subsection{Cost estimates}

The different choices of test and trial functions have an impact also on the computational cost of each method. 
For high order methods, it is convenient to avoid assemblying the full system matrix, and to compute only its action on a vector. This can be achieved by a sequence of matrix-matrix products involving local element matrices.
Calling $Q$ and $n$ the number of quadrature points and the number of basis functions in each direction for each element, we have the following cost entries~\cite{DevilleFischerMund}:
\begin{itemize}
  \item $28Q^2n^4$ for the matrix-matrix multiplications of Equation~\eqref{eq:def_tensor_mat}, for the case of B-spline and NURBS Galerkin Methods;
  \item $20Q^2n^4$ for the matrix-matrix multiplications of Equation~\eqref{eq:def_tensor_mat}, for the Legendre and Chebyshev Galerkin Methods;
  \item $28n^4$ for the matrix-matrix multiplications of Equation~\eqref{eq:LapB_collocation}, for the B-spline and NURBS Collocation Methods;
  \item $16n^4$ for the matrix-matrix multiplications of Equation~\eqref{eq:LapB_collocation}, for the Chebyshev Collocation Method.
\end{itemize}
In most practical codes, the matrix-matrix multiplications are the subroutines absorbing most of the computational resources and time. For this reason, we proceed discussing the computational complexity of this performance-critical phase.

Following the common practice, in the quadrature rules we take $Q=n$ for the Legendre and Chebyshev methods, and $Q=2n-1$ for the B-spline and NURBS methods.

\begin{figure}[tbp]
  \centering
\subfloat[]{
  \includegraphics[width=0.45\textwidth]{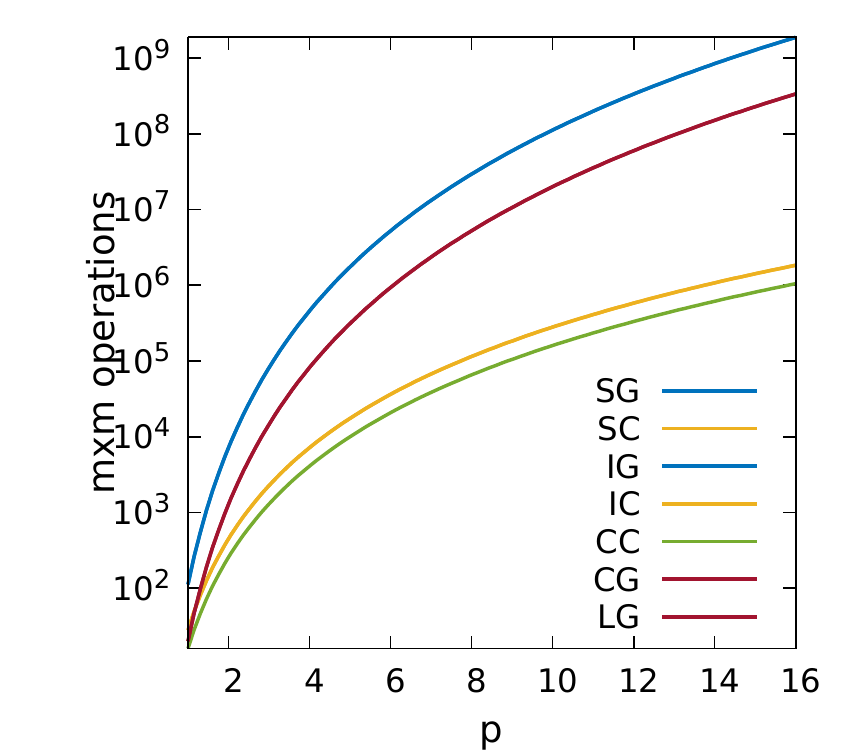}
}
\caption{Operation count for an elementwise matrix-matrix multiplication, for the 7 numerical methods of Section~\ref{sec:basis}. Note that some of the methods share the same color, since the operation count is the same and the relative curves overlap. More precisely, the SG and IG are both shown in blue, SC and IC are both in yellow, CG and LG are both red.}
  \label{fig:operations}
\end{figure}

As a result, for a fixed polynomial order $p$, we have a leading-order operation count for the matrix-matrix multiplication that differs quite significantly for the different numerical methods, that we summarize in Table~\ref{tab:operations}. The same information is visualized in Figure~\ref{fig:operations}.

\begin{table}
  \caption{Leading-order cost for an elementwise matrix-matrix multiplication, for the 7 numerical methods of Section~\ref{sec:basis}.}
  \label{tab:operations}
  \centering
  \begin{tabular}{ccccccc}
    \toprule
    SG      & SC     & IG      & IC    & CC     & CG    & LG \\
    \midrule
    $112n^6$ & $28n^4$ & $112n^6$ & $28n^4$ & $16n^4$ & $20n^6$ & $20n^6$ \\
    \bottomrule
  \end{tabular}
\end{table}

While this choice guarantees exact polynomial integration, there is a vast literature dedicated to better choices of quadrature formulas for Isogeometric methods that shows how the computational cost argument could be less stringent than the one depicted here (see for example \cite{Hughes2010a,Auricchio2012d,Schillinger2014}). We point out that the coefficients in Table~\ref{tab:operations} could be smaller for the B-spline and Isogeometric methods, see in particular the new memory-efficient assembly strategies introduced in~\cite{CALABRO2017606}.

\section{Numerical results}
\label{sec:numerical_results}
In this section we compare the numerical methods presented in Section~\ref{sec:methods} for the Laplace--Beltrami and Allen--Cahn equations on moderately complex surfaces.
The numerical experiments are designed to assess the behaviour of the different methods with respect to order elevation.
Due to the high flexibility in the definition of B-spline and NURBS basis functions, there are many ways to construct and refine such basis.
Two common ways to increase the order of B-spline and NURBS functions are $p$ and $k$-refinement (see e.g.~\cite{IGABook},~\cite{IGAFirst},~\cite{Sangallihpk}). While for $p$-refinement there is a definition the literature agrees on, namely the increase by one of all the knots' multiplicity, a single step of $k$-refinement may consist of a $p$-refinement followed by some (somewhat arbitrary) knot insertions, at locations that do not coincide with existing knots.
The definition of a single $k$-refinement step adopted here consists in increasing by one the multiplicity of both the internal and the end knots, followed by the insertion of one internal knot per knot interval, with multiplicity equal to one. Subsequent $k$-refinements do not start from previous stages, but from the knot vector of the original geometry, i.e., the internally inserted knots are removed before elevating the degree of the B-splines.

\begin{figure}[tbp]
  \centering
  \begin{minipage}{.31\textwidth}
  \begin{tikzpicture}[scale=2]
    \draw [thick] (0,0) -- (2,0);
    \draw [fill] (0,0) circle [radius=0.05];
    \draw [fill] (0,0.15) circle [radius=0.05];
    \draw [fill] (0,0.30) circle [radius=0.05];
    \draw [fill] (1,0) circle [radius=0.05];
    \draw [fill] (1,0.15) circle [radius=0.05];
    \draw [fill] (2,0) circle [radius=0.05];
    \draw [fill] (2,0.15) circle [radius=0.05];
    \draw [fill] (2,0.30) circle [radius=0.05];
  \end{tikzpicture}
  \end{minipage}
  \begin{minipage}{.31\textwidth}
\begin{tikzpicture}[scale=2]
    \draw [thick] (0,0) -- (2,0);
    \draw [fill] (0,0) circle [radius=0.05];
    \draw [fill] (0,0.15) circle [radius=0.05];
    \draw [fill] (0,0.30) circle [radius=0.05];
    \draw [fill] (0,0.45) circle [radius=0.05];
    \draw [fill] (0.5,0) circle [radius=0.05];
    \draw [fill] (1,0) circle [radius=0.05];
    \draw [fill] (1,0.15) circle [radius=0.05];
    \draw [fill] (1,0.30) circle [radius=0.05];
    \draw [fill] (1.5,0) circle [radius=0.05];
    \draw [fill] (2,0) circle [radius=0.05];
    \draw [fill] (2,0.15) circle [radius=0.05];
    \draw [fill] (2,0.30) circle [radius=0.05];
    \draw [fill] (2,0.45) circle [radius=0.05];
  \end{tikzpicture}
  \end{minipage}
\begin{minipage}{.31\textwidth}
\begin{tikzpicture}[scale=2]
    \draw [thick] (0,0) -- (2,0);
    \draw [fill] (0,0) circle [radius=0.05];
    \draw [fill] (0,0.15) circle [radius=0.05];
    \draw [fill] (0,0.30) circle [radius=0.05];
    \draw [fill] (0,0.45) circle [radius=0.05];
    \draw [fill] (0,0.60) circle [radius=0.05];
    \draw [fill] (0.333,0) circle [radius=0.05];
    \draw [fill] (0.667,0) circle [radius=0.05];
    \draw [fill] (1,0) circle [radius=0.05];
    \draw [fill] (1,0.15) circle [radius=0.05];
    \draw [fill] (1,0.30) circle [radius=0.05];
    \draw [fill] (1,0.45) circle [radius=0.05];
    \draw [fill] (1.333,0) circle [radius=0.05];
    \draw [fill] (1.667,0) circle [radius=0.05];
    \draw [fill] (2,0) circle [radius=0.05];
    \draw [fill] (2,0.15) circle [radius=0.05];
    \draw [fill] (2,0.30) circle [radius=0.05];
    \draw [fill] (2,0.45) circle [radius=0.05];
    \draw [fill] (2,0.60) circle [radius=0.05];
  \end{tikzpicture}
  \end{minipage}
  \caption{Evolution of a knot vector under $k$-refinement. Each dot represents a knot, with the convention that a vertical array of $m$ dots stands for a knot of multiplicity $m$. The initial knot vector (left) is given. The knot vector in the middle is obtained after $k$-refinement of the initial knot vector. The knot vector on the right is obtained after two levels of $k$-refinement.} 
  \label{fig:k_refinement}
\end{figure}
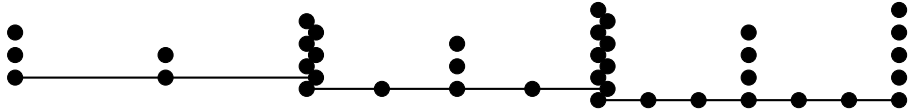

In particular,  given an initial knot vector defining the geometry, in this work we construct the new knot vector obtained after $m$ steps of $k$-refinements by taking the same knot vector, with the multiplicities of every knot increased $m$ times, union with a vector of $m$ new knots of multiplicity one for each knot interval, located in new points equally spaced between each couple of subsequent knots in the original vector. An example that shows our convention is available in Figure~\ref{fig:k_refinement}.

We remark that $k$-refinement can be interpreted as a composition of a $p$-refinement obtained by augmenting the geometry knots' multiplicity, and a series of knot insertions, obtained by inserting new inner knots with multiplicity one.

Both $p$ and $k$-refinement are a form of degree elevation, with the important difference that during $p$-refinement the basis functions keep their original global continuity, while after $m$ steps of $k$-refinement, new internal knots are inserted in the knot vector, generating basis functions at the inserted knots that have the maximum available continuity. For example, if initially the basis functions have order $q$ and are globally $\CC^s$, after $m$ steps of $p$-refinement the new basis functions will have order $q+m$, and maintain the global regularity $\CC^s$. The same initial basis functions, after $m$ steps of $k$-refinement will have order $q+m$, global regularity $\CC^s$, and local regularity $\CC^{q+m-1}$ on the newly inserted knots.

An alternative approach to the $k$-refinement strategy described above, consists in  removing internal knots while increasing the degree of the B-splines. This strategy, also referred to as $k$-coarsening, implies an increase in the global regularity of the B-spline basis functions, at the price of generating non-nested spaces. This approach has the advantage that the number of degrees of freedom does not grow too fast with the polynomial degree, and may also lead to better conditioned matrices than those obtained by keeping fixed the global continuity of the basis functions. A major disadvantage of this approach is related to the fact that the geometry cannot be preserved through $k$-coarsening, making it only useful for trivial geometries, and requiring an additional geometry reconstruction step, which may not be well posed, or may give unsatisfactory results. In this work we only show the $k$-refinement strategy illustrated above, which is guaranteed to preserve the exact geometry.


In all cases the linear systems are assembled with the \texttt{numpy}~\cite{numpy} and \texttt{igakit}~\cite{igakit} libraries. 

\subsection{Description of the test cases}

A comprehensive test requires that the following possibilities are fully considered:
\begin{itemize}
  \item flat surfaces (domains) and curved surfaces;
  \item collocation and Galerkin;
  \item mixed boundary conditions;
  \item $p$-refinement and $k$-refinement (when applicable);
  \item linear and nonlinear problems.
\end{itemize}

In addition, it is sensible to check that for B-spline and NURBS collocation methods, the convergence rate does not depend drastically on the choice of collocation points.

The test cases we set up consist in a homogeneous Laplace--Beltrami problem on two geometries of increasing complexity, including a comparison between different collocation strategies, the analysis of a mixed Neumann--Dirichlet boundary condition problem, and a nonlinear test case solving the Allen--Cahn equation.

In addition to standard $p$-refinement strategies, we also make a comparison between hybrid NURBS-SEM methods and a variant of the isogeometric $k$-refinement strategy. For consistency, in this case the comparison is not performed on the basis of the polynomial degree, but on the total number of degrees of freedom, since the $k$-refinement strategy that we employ introduces a higher number of degrees of freedom for a fixed polynomial order than standard $p$-refinement in SEM.

The first surface we consider consists of a domain of $\mathbb{R}^2$, lying on a plane oblique to the three coordinate axis of $\mathbb{R}^3$. Since the surface is flat, any harmonic function in the plane of the surface will also solve the Laplace--Beltrami equation.
We consider a transcendental harmonic function as a reference solution, whose expression is reported in Table~\ref{tab:experiments}.
The restriction of this function to the surface boundary provides the required Dirichlet data.

The flat geometry that we consider is a quarter of annulus, shown in Figure~\ref{fig:quarter_annulus}, obtained by considering the region between two concentric circles with inner radius $R_1=0.5$ and outer radius $R_2=1$, and two orthogonal diameters.

The annular surface is represented by a mesh of $2\times2$ elements, as shown in Figure~\ref{fig:quarter_elements} (a). The collocation points for the Chebyshev (GLC nodes) and Isogeometric (Greville nodes) methods of degree 7 are shown in Figures~\ref{fig:quarter_elements} (b) and (c) respectively.

\begin{figure}[tbp]
  \centering
  \begin{tikzpicture}[scale=2.5]
    \draw[thick] (0.5,0) -- (1,0);
    \draw[thick] (0,0.5) -- (0,1);
    \draw[thick] (0.5,0) arc (0:90:0.5);
    \draw[thick] (1,0) arc (0:90:1);
    \draw[<->] (0,-0.3) -- (0.5,-0.3) node[midway,above]{$R_1$};
    \draw[<->] (0,-0.6) -- (1,-0.6) node[midway,above]{$R_2$};
    \draw (0.5,0) -- (0.5,-0.3);
    \draw (0,0) -- (0,-0.6);
    \draw (1,0) -- (1,-0.6);
    \draw (0,-0.05) -- (0,0.05);
    \draw (-0.05,0) -- (0.05,0);
  \end{tikzpicture}
  \caption{Surface geometry for the test case n.1.}
  \label{fig:quarter_annulus}
\end{figure}
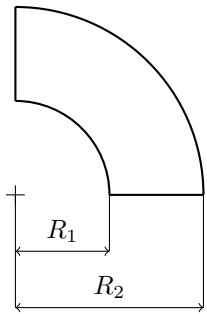
\begin{figure}[tbp]
  \centering
\subfloat[]{
  \begin{tikzpicture}[scale=4.] 
    \draw[thick] (0.5,0) -- (1,0);
    \draw (0.3536,0.3536) -- (0.7071,0.7071);
    \draw[thick] (0,0.5) -- (0,1);
    \draw[thick] (0.5,0) arc (0:90:0.5);
    \draw        (0.75,0) arc (0:90:0.75);
    \draw[thick] (1,0) arc (0:90:1);
  \end{tikzpicture}
}\qquad
\subfloat[]{
   \begin{tikzpicture}[scale=4.]
    \draw[thick] (0.5,0) -- (1,0);
    \draw (0.3536,0.3536) -- (0.7071,0.7071);
    \draw[thick] (0,0.5) -- (0,1);
    \draw[thick] (0.5,0) arc (0:90:0.5);
    \draw        (0.75,0) arc (0:90:0.75);
    \draw[thick] (1,0) arc (0:90:1);
    \draw        (0.512379,0) arc (0:90:0.512379);
    \draw        (0.547064,0) arc (0:90:0.547064);
    \draw        (0.597185,0) arc (0:90:0.597185);
    \draw        (0.652815,0) arc (0:90:0.652815);
    \draw        (0.702936,0) arc (0:90:0.702936);
    \draw        (0.737621,0) arc (0:90:0.737621);
    \draw        (0.762379,0) arc (0:90:0.762379);
    \draw        (0.797064,0) arc (0:90:0.797064);
    \draw        (0.847185,0) arc (0:90:0.847185);
    \draw        (0.902815,0) arc (0:90:0.902815);
    \draw        (0.952936,0) arc (0:90:0.952936);
    \draw        (0.987621,0) arc (0:90:0.987621);
    \draw (0.499622, 0.0194398) -- (0.999244, 0.0388796);
    \draw (0.494545, 0.0736585) -- (0.989089, 0.147317);
    \draw (0.476876, 0.150297) -- (0.953752, 0.300594);
    \draw (0.443478, 0.230926) -- (0.886957, 0.461853);
    \draw (0.40178, 0.297611) -- (0.803561, 0.595223);
    \draw (0.367032, 0.33954) -- (0.734064, 0.67908);
    \draw (0.0194398, 0.499622) -- (0.0388796, 0.999244);
    \draw (0.0736585, 0.494545) -- (0.147317, 0.989089);
    \draw (0.150297, 0.476876) -- (0.300594, 0.953752);
    \draw (0.230926, 0.443478) -- (0.461853, 0.886957);
    \draw (0.297611, 0.40178) -- (0.595223, 0.803561);
    \draw (0.33954, 0.367032) -- (0.67908, 0.734064);
  \end{tikzpicture}
}\qquad
\subfloat[]{
  \begin{tikzpicture}[scale=4.]
    \draw[thick] (0.5,0) -- (1,0);
    \draw[thick] (0,0.5) -- (0,1);
    \draw[thick] (0.5,0) arc (0:90:0.5);
    \draw[thick] (1,0) arc (0:90:1);
    \draw (0.53125,0) arc (0:90:0.53125);
    \draw (0.56250,0) arc (0:90:0.56250);
    \draw (0.59375,0) arc (0:90:0.59375);
    \draw (0.62500,0) arc (0:90:0.62500);
    \draw (0.65625,0) arc (0:90:0.65625);
    \draw (0.68750,0) arc (0:90:0.68750);
    \draw (0.71875,0) arc (0:90:0.71875);
    \draw (0.75000,0) arc (0:90:0.75000);
    \draw (0.78125,0) arc (0:90:0.78125);
    \draw (0.81250,0) arc (0:90:0.81250);
    \draw (0.84375,0) arc (0:90:0.84375);
    \draw (0.87500,0) arc (0:90:0.87500);
    \draw (0.90625,0) arc (0:90:0.90625);
    \draw (0.93750,0) arc (0:90:0.93750);
    \draw (0.96875,0) arc (0:90:0.96875);
    \draw (0.497592, 0.0490086) -- (0.995185, 0.0980171);
    \draw (0.490393, 0.0975452) -- (0.980785, 0.19509);
    \draw (0.47847, 0.145142) -- (0.95694, 0.290285);
    \draw (0.46194, 0.191342) -- (0.92388, 0.382683);
    \draw (0.440961, 0.235698) -- (0.881921, 0.471397);
    \draw (0.415735, 0.277785) -- (0.83147, 0.55557);
    \draw (0.386505, 0.317197) -- (0.77301, 0.634393);
    \draw (0.353553, 0.353553) -- (0.707107, 0.707107);
    \draw (0.317197, 0.386505) -- (0.634393, 0.77301);
    \draw (0.277785, 0.415735) -- (0.55557, 0.83147);
    \draw (0.235698, 0.440961) -- (0.471397, 0.881921);
    \draw (0.191342, 0.46194) -- (0.382683,0.92388);
    \draw (0.145142, 0.47847) -- (0.290285,0.95694);
    \draw (0.0975452, 0.490393) -- (0.195059,0.980785);
    \draw (0.0490086, 0.497592) -- (0.0980171,0.995185);
  \end{tikzpicture}
}
  \caption{Discretization of the quarter of annulus in four elements (a); Gauss--Lobatto--Chebyshev collocation points (b); Greville collocation points (c).}
  \label{fig:quarter_elements}
\end{figure}
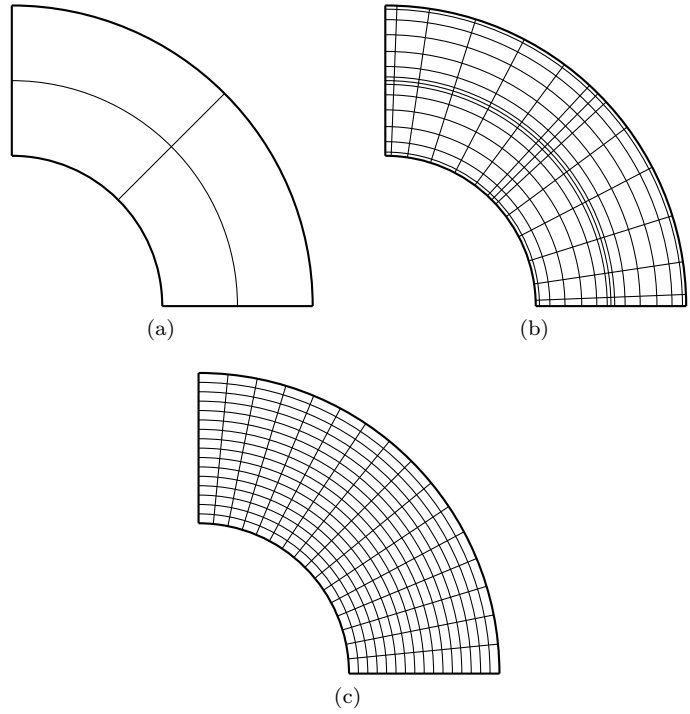

\begin{figure}[tbp]
  \centering
\subfloat[]{
  \includegraphics[width=0.35\textwidth]{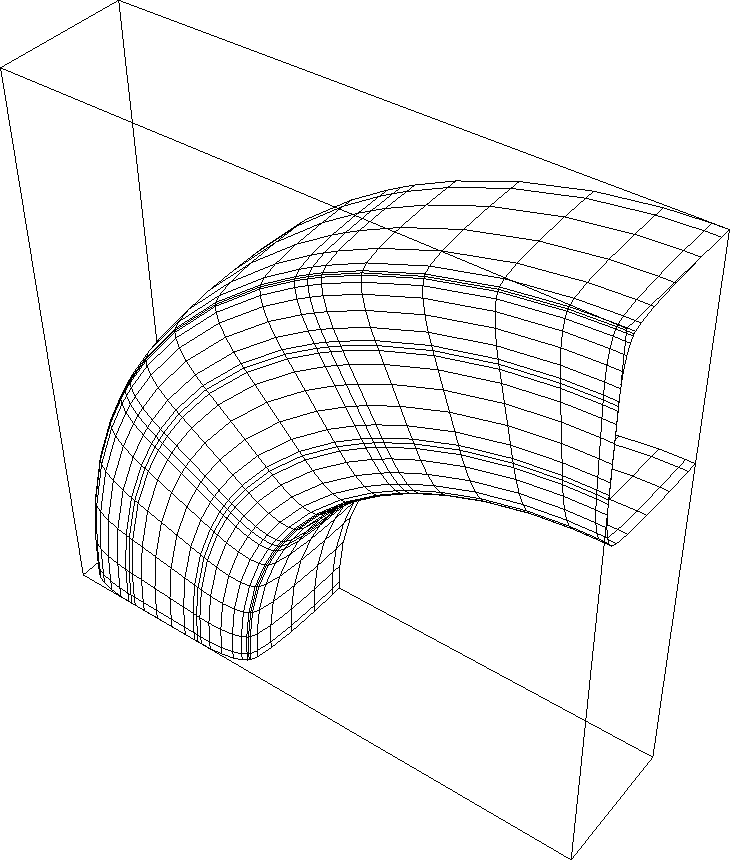}
}\qquad
\subfloat[]{
  \includegraphics[width=0.35\textwidth]{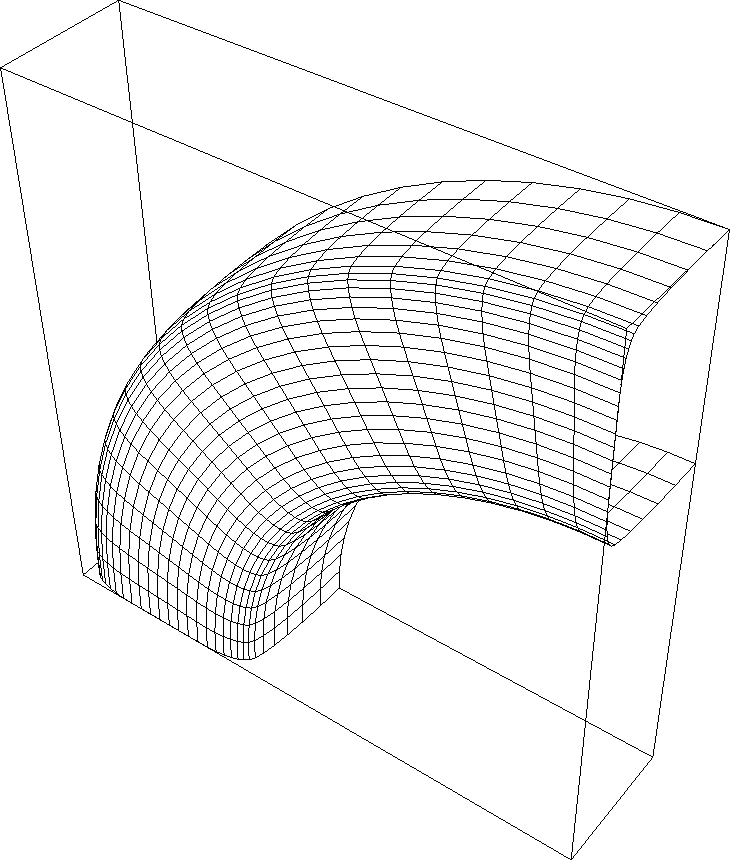}
}
  \caption{Discretization in $3\times5$ elements of the C-shaped surface for the second test case: Gauss--Lobatto--Chebyshev collocation points (a); Greville collocation points (b).}
  \label{fig:rotated_elements}
\end{figure}

The curved geometry considered for this test is shown in Figure~\ref{fig:rotated_elements}, and is obtained by revolving a C-shaped profile around a quarter of a circle of diameter equal to 5 times the sectional heigth.
Since this surface has nonzero curvature, an harmonic function in $\mathbb{R}^n$ will not in general be a solution of the Laplace--Beltrami equation.
Consequently, in this case the numerical results are compared against a manufactured solution, as reported in Table~\ref{tab:experiments}.
For this geometry we construct a mesh with 3 elements in the direction of revolution, and in 5 elements in the radial direction.
In Figure~\ref{fig:rotated_elements} we show the collocation points for a Chebyshev collocation method (a) and for a B-spline or NURBS collocation method (b), both of degree 7.

The full set of numerical experiments carried out in the following sections is summarized in Table~\ref{tab:experiments}.
\begin{table}
  \caption{Summary of the numerical experiments performed in this work to compare the high-order methods of Section~\ref{sec:methods}.}
  \label{tab:experiments}
  \centering
  \begin{tabular}{clcl}
    \toprule
    $n$ & description & geometry & solution $u(\xx)$ \\
    \midrule
     1 & $p$-ref flat & Annulus & $-\frac{1}{2\pi}\log|\xx-(1,1,0)|$ \\
     2 & collocation iga & Annulus &  $-\frac{1}{2\pi}\log|\xx-(1,1,0)|$ \\
     3 & $p$-ref mixed bc & Annulus &  $-\frac{1}{2\pi}\log|\xx-(1,1,0)|$ \\ 
     4 &  $p$-ref curved &C-surface &  $\cos(x_2)\cos(x_3)$ \\
    5& {$k$-ref flat} & {Annulus} & $-\frac{1}{2\pi}\log|\xx-(1,1,0)|$ \\
    6 & $k$-ref nonlinear & {Annulus} & $x_1^2-x_2^3$ \\
    \bottomrule
  \end{tabular}
\end{table}

\subsection{$p$-refinement, flat geometry}
\label{sec:quarter_annulus}
\begin{figure}[tbp]
  \centering
  \subfloat[]{
    \includegraphics[width=0.45\textwidth]{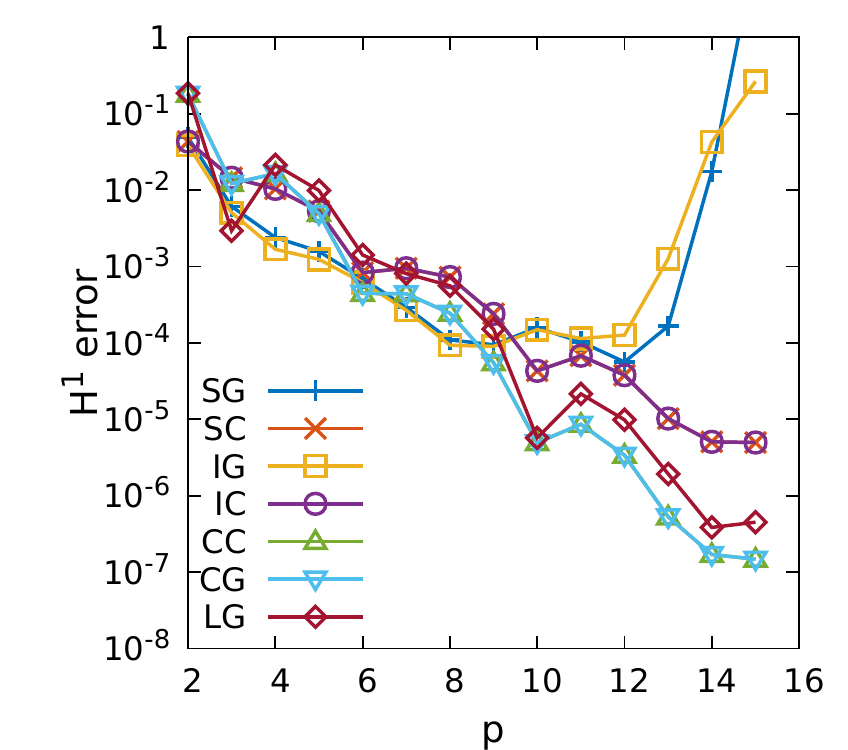}
  } 
  \subfloat[]{
    \includegraphics[width=0.45\textwidth]{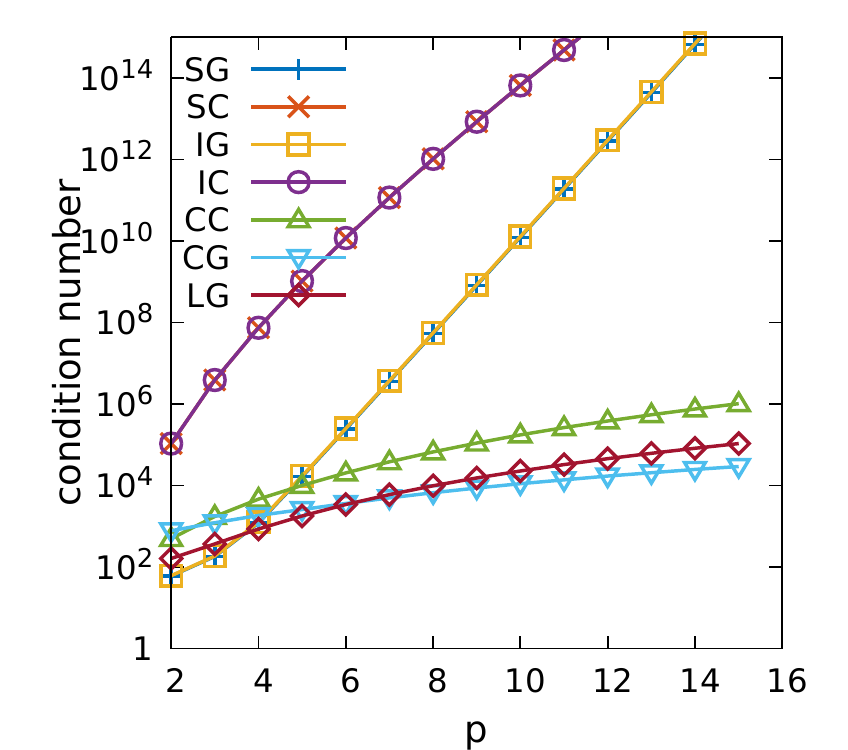}
  }
  \caption{Comparison of the numerical methods for the annular geometry and the transcendental solution $u(\xx)=\log|\xx-(1,1,0)|$. Error in the $H^1$ norm (left) and matrix condition number (right) as function of the polynomial degree.}
  \label{fig:results_annulus_log}
\end{figure}

The results for the case with annular geometry and a transcendental harmonic solution are shown in Figure~\ref{fig:results_annulus_log}. 
Figure~\ref{fig:results_annulus_log} (a) shows that for a polynomial degree between 2 and 9 all the methods considered here achieve the expected exponential convergence. For polynomial degrees up to 11, the B-spline and NURBS method except for the IGA-Collocation method fail at keeping the exponential trend, and for even higher polynomial degrees, only the Spectral Element methods show a satisfactory behaviour, reaching spectral accuracy.

The observed behaviour can be explained by looking at the matrix condition number as a function of the polynomial degree $p$, shown in Figure~\ref{fig:results_annulus_log} (b). From this picture, it is clear that B-spline and NURBS matrices are too ill conditioned for high orders. Conversely, the condition number of the Chebyshev and Legendre matrices increases only algebraically, with the expected asymptotic order $O(p^4)$.

\subsection{Dependence on the choice of collocation points}
The problem set up with annular geometry and transcendental solution is at the basis of two further tests regarding respectively the choice of collocation points in Isogeometric analysis and the use of mixed Neumann--Dirichlet boundary conditions.

\begin{figure}[tbp]
  \centering
  \subfloat[]{
    \includegraphics[width=0.45\textwidth]{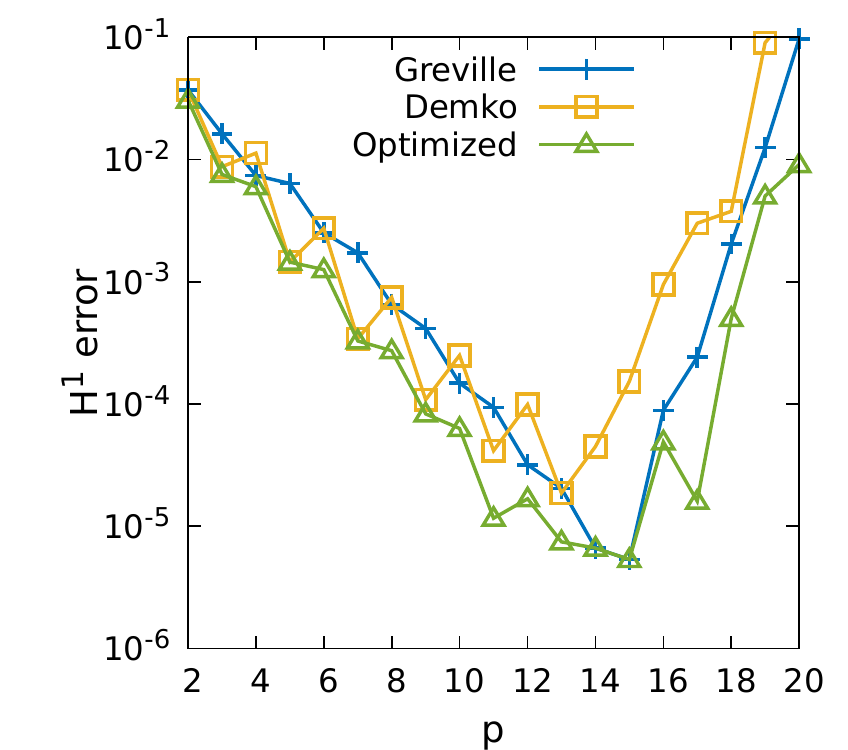}
  }
  \subfloat[]{
    \includegraphics[width=0.45\textwidth]{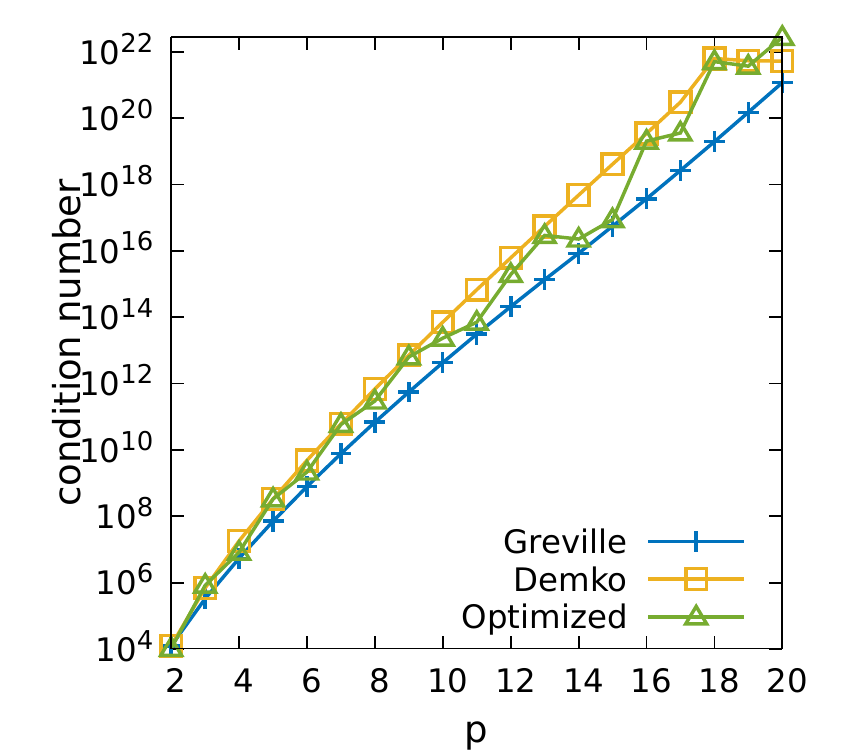}
  }
  \caption{Comparison of the three Isogeometric collocation methods based on three sets of collocation points. Error in the $H^1$ norm (left) and matrix condition number (right) as function of the polynomial degree.}
  \label{fig:collocation}
\end{figure}

To this end, we repeat the test described above for the Isogeometric collocation method with three different sets of collocation points, namely:
\begin{description}
  \item[Greville points] Greville points are defined as the points in the unit interval where each basis function achieves its local maximum.
  \item[Demko points] also called \emph{Demko abscissae}, were introduced in~\cite{Demko} and applied in~\cite{Auricchio2012} to Isogeometric Collocation. These are the points at which Chebyshev splines (i.e. splines which oscillate between $-1$ and $1$) achieve a maximum or a minimum. By definition, the maxima and minima of Chebyshev splines are exactly $1$ and $-1$ respectively.
  \item[Optimized points] In this case, a minimizer package from the library \texttt{scipy}~\cite{scipy} is used to find the collocation points that minimize the $H^1$ norm of the error for each degree $p$ of the basis functions. This set of points clearly is out of reach if the solution is not known in advance, or if a sharp error estimator is not available. However, it is instructive to compute anyways this set of points as a lower bound for the comparison with Greville and Demko points.
\end{description}
Only for this test, due to the very high computational cost of the optimization subroutine, we subdivide the geometry in a single large element instead of four smaller elements.
The results of this comparison, shown in Figure~\ref{fig:collocation}, show that on average Demko points may lead to slightly smaller errors than Greville points, but the convergence rate and the maximum attainable accuracy seem not to be much influenced by the choice of collocation points. However, it is worth mentioning that many other choices of collocation points are available in the literature (see, for example,~\cite{anitescu2015isogeometric}), and that the list we provide here is by no means complete.

\subsection{Problems with mixed boundary conditions}
To check the treatment of Neumann boundary conditions, we replace the Dirichlet boundary condition in the curved edge of radius $R_2=1$ with Neumann boundary conditions.

\begin{figure}[tbp]
  \centering
  \subfloat[]{
    \includegraphics[width=0.45\textwidth]{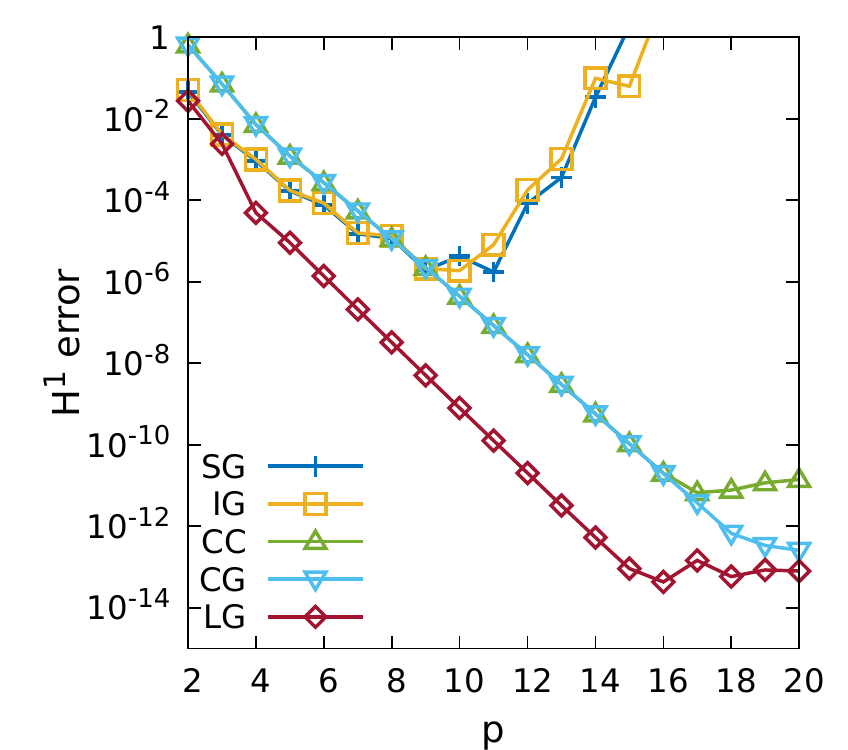}
  }
  \subfloat[]{
    \includegraphics[width=0.45\textwidth]{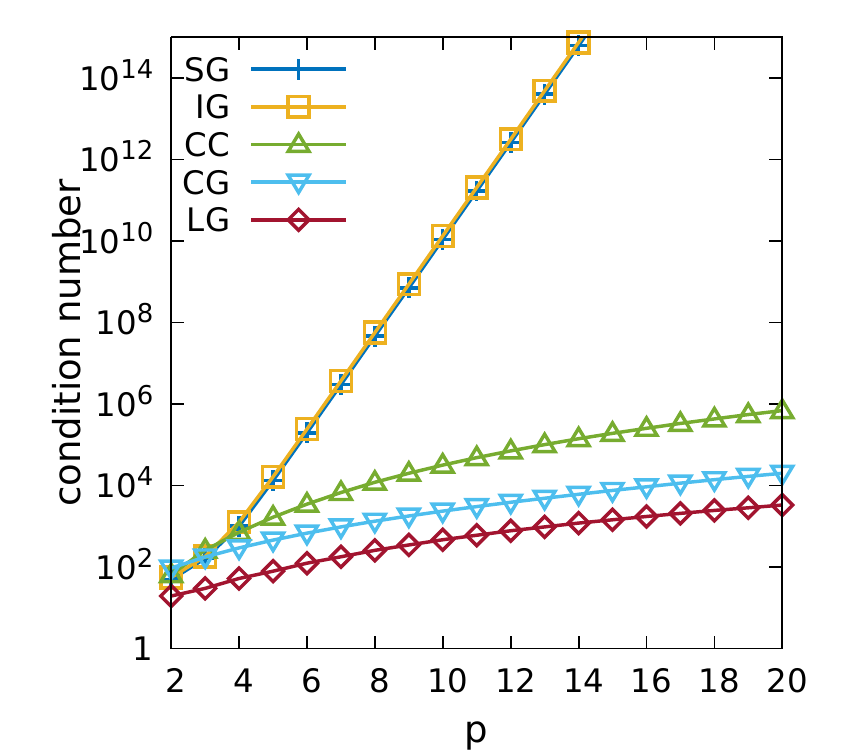}
  }
  \caption{Comparison of the numerical methods for the annular geometry with transcendental solution and Neumann boundary conditions.  Error in the $H^1$ norm (left) and matrix condition number (right) as function of the polynomial degree.}
  \label{fig:annulus_neumann}
\end{figure}

The results of this last test case, shown in Figure~\ref{fig:annulus_neumann}, confirm the good behaviour of Spectral Element methods with respect to degree elevation. The imposition of natural boundary conditions in Galerkin methods appears to be better conditioned than in collocation methods.

\subsection{$p$-refinement, curved geometry}
\label{sec:revolved_c_beam}

\begin{figure}[tbp]
  \centering
  \subfloat[]{
    \includegraphics[width=0.45\textwidth]{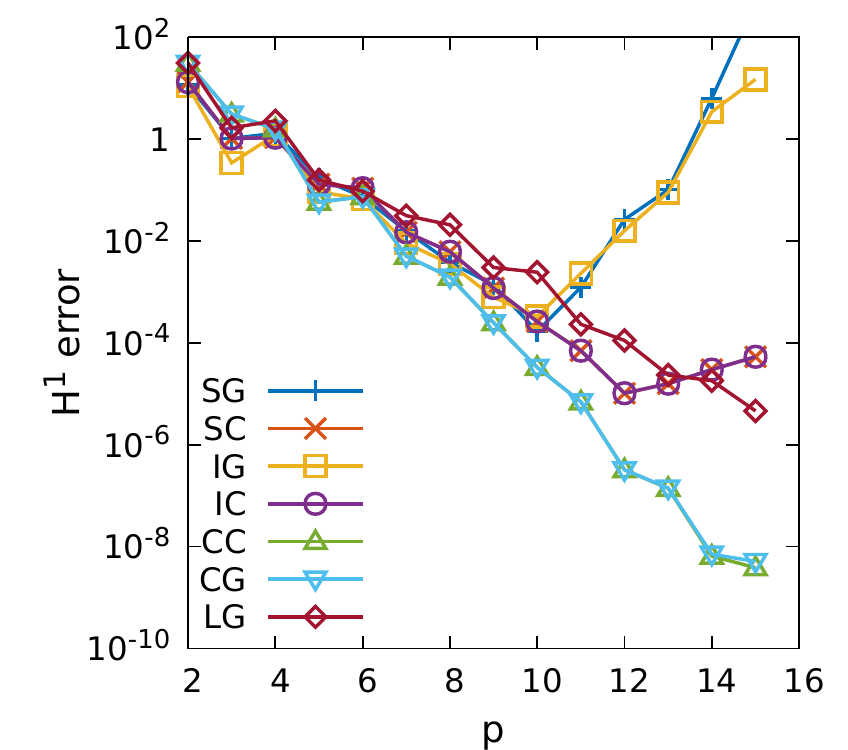}
  } 
  \subfloat[]{
     \includegraphics[width=0.45\textwidth]{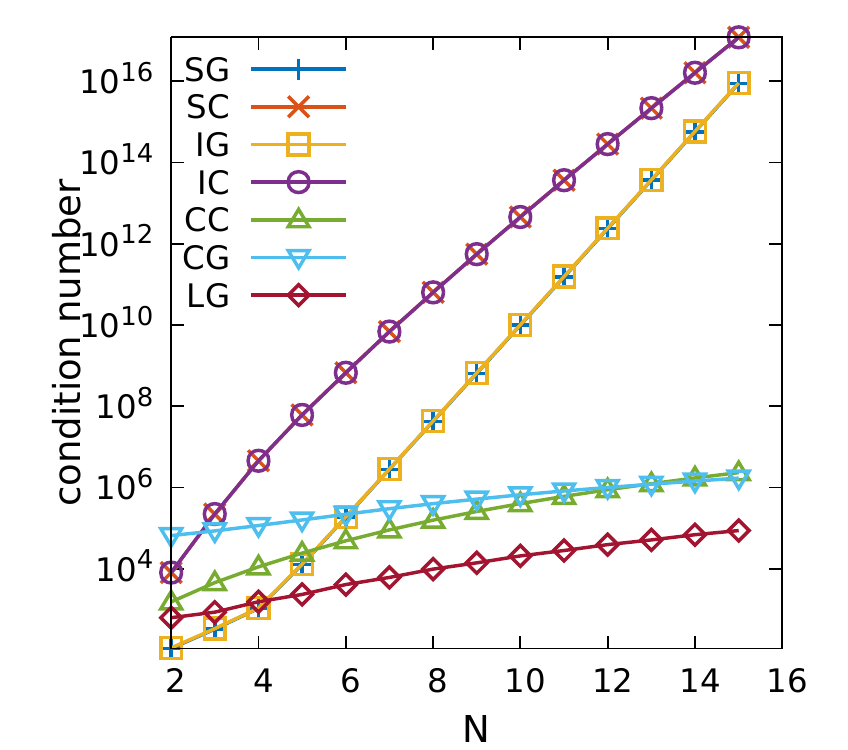}
   }
  \caption{Comparison of the numerical methods for the C-shaped surface and the transcendental solution $u(\xx)=\cos(x_2)\cos(x_3)$.  Error in the $H^1$ norm (left) and matrix condition number (right) as function of the polynomial degree.}
  \label{fig:results_cbeam_cos}
\end{figure}
\begin{figure}[htb]
  \centering
  \subfloat[]{
    \includegraphics[width=0.48\textwidth]{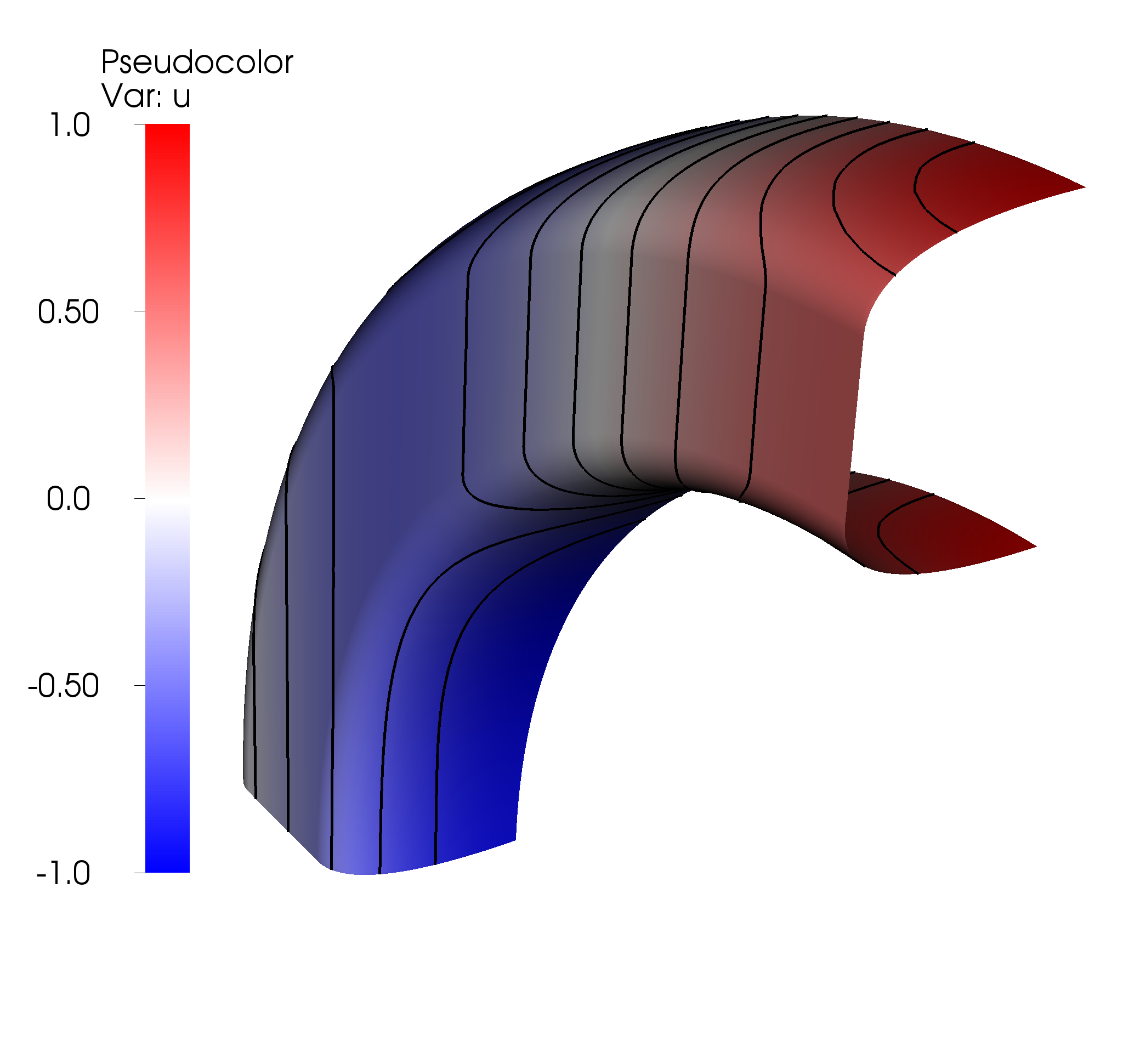}
  }
  \subfloat[]{
    \includegraphics[width=0.48\textwidth]{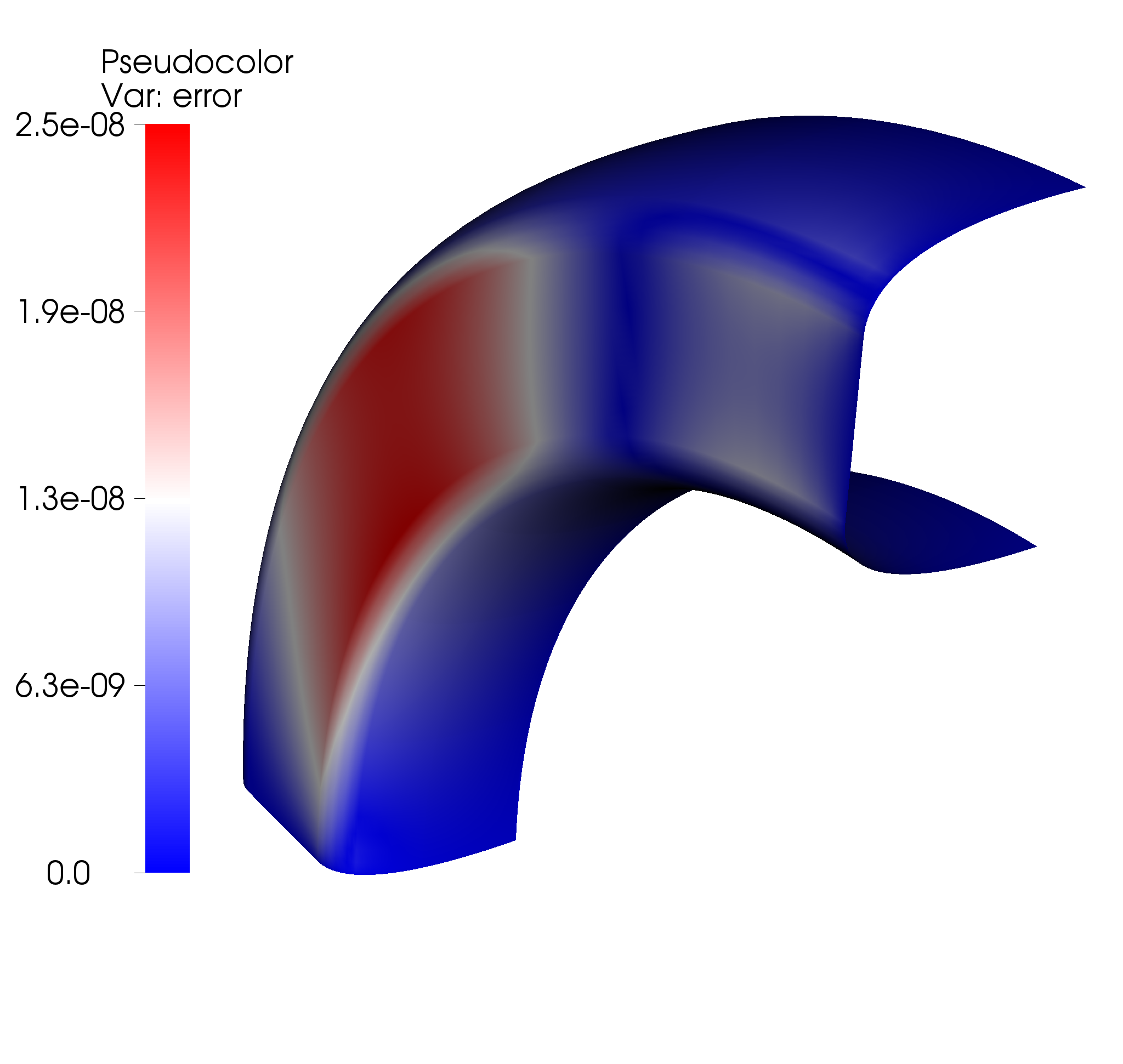}
  }
  \caption{Plot of the solution $u(\xx)=\cos(x_2)\cos(x_3)$ (a) and of the pointwise error for the LG method with polynomial degree equal to $15$ (b).}
  \label{fig:testl}
\end{figure}

The results using the manufactured  solution given in Table~\ref{tab:experiments} are shown in Figure~\ref{fig:results_cbeam_cos}, where we plot on the left the $H^1$ error and on the right the matrix condition number as a function of the polynomial degree. We report in Figure~\ref{fig:testl} a plot of the solution (a) together with a map of the pointwise error (b) for the Legendre case of order 15.
As for the previous experiments, we see that B-spline and NURBS methods fail to achieve spectral accuracy, while performing very well at relatively low polynomial degrees.

\subsection{$k$-refinement, flat geometry}
\label{sec:comparison_k_and_p}
In view of the results shown in Sections~\ref{sec:quarter_annulus} and~\ref{sec:revolved_c_beam}, it may be tempting to conclude that IGA underperforms in terms of $p$-refinement if compared with a hybrid NURBS-mapped Spectral Element Method.
However, any comparison of numerical methods focusing on the behaviour with respect to degree elevation of the underlying basis functions would not be complete if $k$-refinement were not considered.
An important feature of our $k$-refinement strategy is that it allows the addition of internal knots while increasing the polynomial degree. This cannot be achieved by simple $p$-refinement, and may give an edge to $k$-refinement over $p$-refinement in terms of accuracy, thanks to these additional degrees of freedom.

The goal of the present section is therefore to compare the performance of increasing the order of the method by performing $m$ steps of $k$-refinement in B-spline and NURBS methods with $m$ standard steps of $p$-refinement in Spectral Element Methods. 

The tests are performed on the quarter of annulus geometry, with initial knot vectors given by:
\begin{equation}
  \Theta^1 = \{ 0,0,0,0.5,0.5,1,1,1\} \qquad
  \Theta^2 = \{ 0,0,0, 1,1,1\}.
  \label{eq:theta_1_theta_2}
\end{equation}
As for the previous cases, we consider a transcendental solution reported in Table~\ref{tab:experiments}. At step $m$, we generate the new knot vectors by increasing the multiplicity of each of the above knots by $m$, elevating the degree of the polynomial base by $m$, and  inserting $2m$ knots in $\Theta^1$, of which $m$ equally spaced between $(0,0.5)$, $m$ equally spaced between $(0.5, 1)$, and finally $m$ new knots in $\Theta^2$,  equally spaced between $(0,1)$.

In all the test cases, Isogeometric methods with $k$-refinement are remarkably efficient for polynomial degrees up to $10$, where on average these deliver a solution two orders of magnitude more accurate than their SEM counterparts with the same polynomial degree. This is to be expected, since the $k$-refinement procedure generates a set of basis functions which is $2m^2$ bigger than the corresponding $p$-refinement SEM basis functions, and makes the comparison between the two methods unfair if done in terms of the polynomial degree alone.

\begin{figure}[tbp]
  \centering
  \subfloat[]{
    \includegraphics[width=0.45\textwidth]{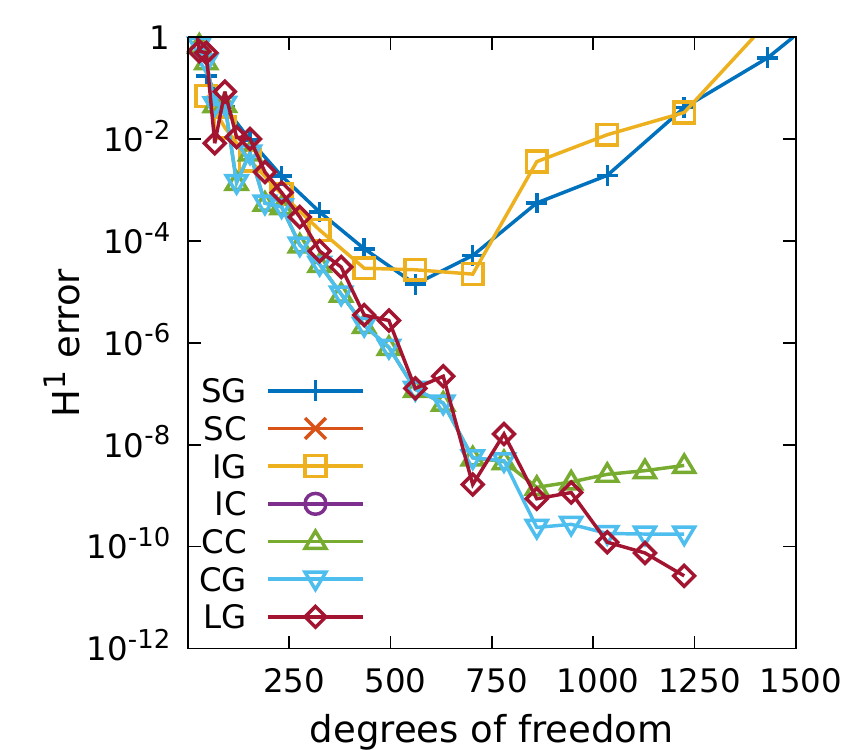}
  }
  \subfloat[]{
    \includegraphics[width=0.45\textwidth]{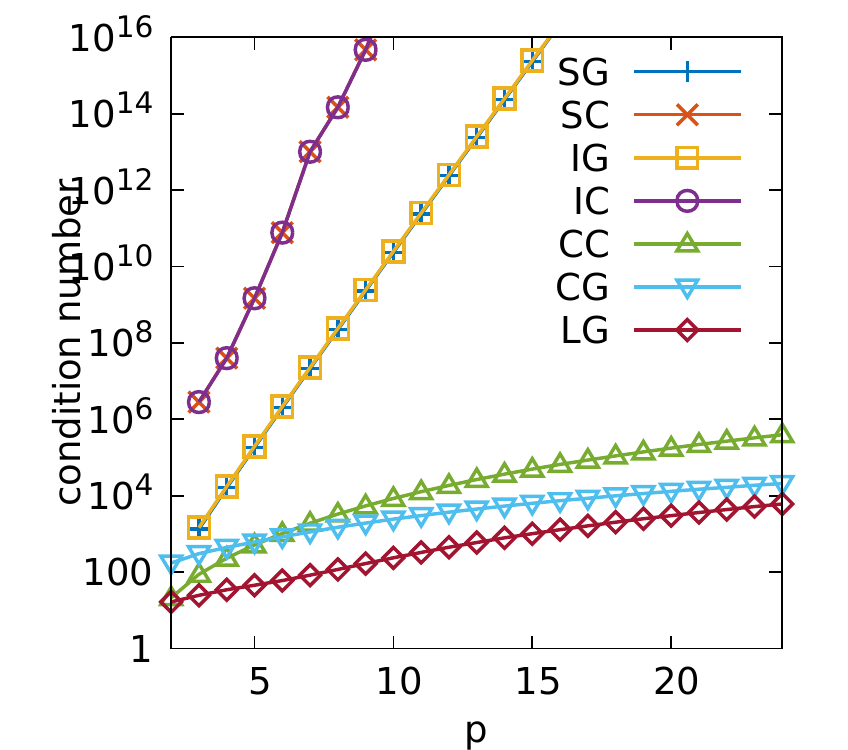}
  }
  \caption{Comparison of the $k$-refinement methods for the annular geometry and the transcendental solution $u(\xx) = \log|\xx - (1, 1, 0)|$.
    Error in the $H^1$ norm as function of the number of degrees of freedom (left) and matrix condition number as function of the polynomial degree (right).}
  \label{fig:Kbeam}
\end{figure}

A more fair comparison is obtained when the error and the condition number are plotted in terms of the number of basis functions, as in Figure~\ref{fig:Kbeam}. In this case, the accuracy per degree of freedom is substantially the same, up to order 10. As the polynomial degree is increased, the condition number of the Isogeometric matrices becomes so high (as shown in Figure~\ref{fig:k_refinement_cond} (a)) that no meaningful solution is delivered already at $p=12$. Conversely, SEM achieved spectral accuracy in all test cases.

Arguably, one could consider a comparison between $hp$-refinement on the spectral methods and $k$-refinement on the isogeometric methods. It is in principle possible to add $m$ internal knots before performing $m$ degree elevation, and to compare this strategy to our $k$-refinement strategy (perform $m$ degree elevation and then add $m$ internal nots). 

In this case, the number of degrees of freedom in the $hp$-refinement would grow much more quickly than our $k$-refinement strategy, and a comparison on the basis of the number of degrees of freedom would lead to a large imbalance of the polynomial degrees for the same number of degrees of freedom, making this comparison less significant. 

\subsection{$k$-refinement, nonlinear problem}

\begin{figure}[tbp]
  \centering
  \subfloat[]{
    \includegraphics[width=0.45\textwidth]{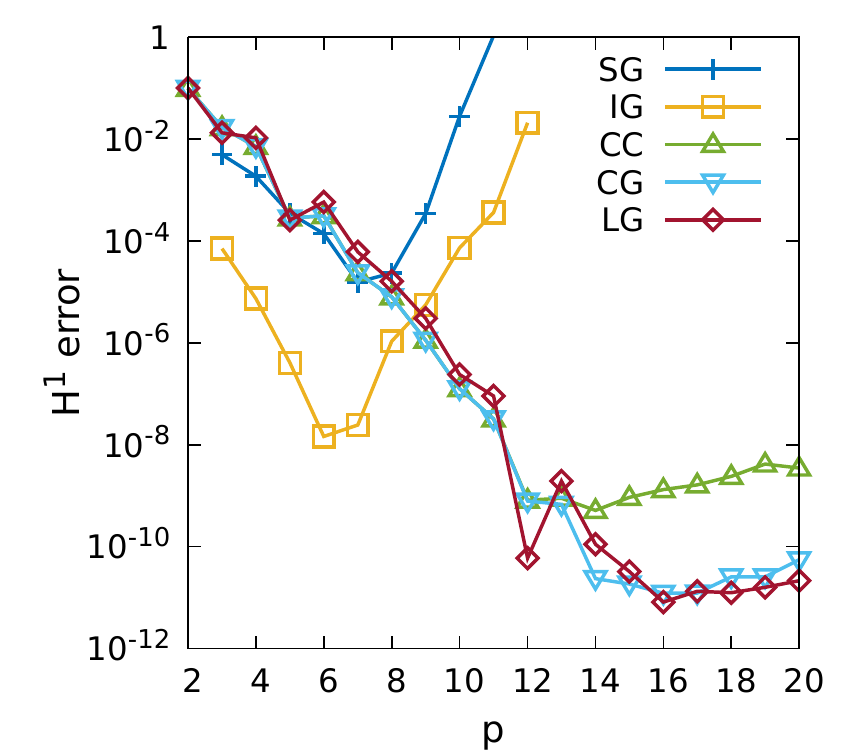}
  }
  \subfloat[]{
    \includegraphics[width=0.45\textwidth]{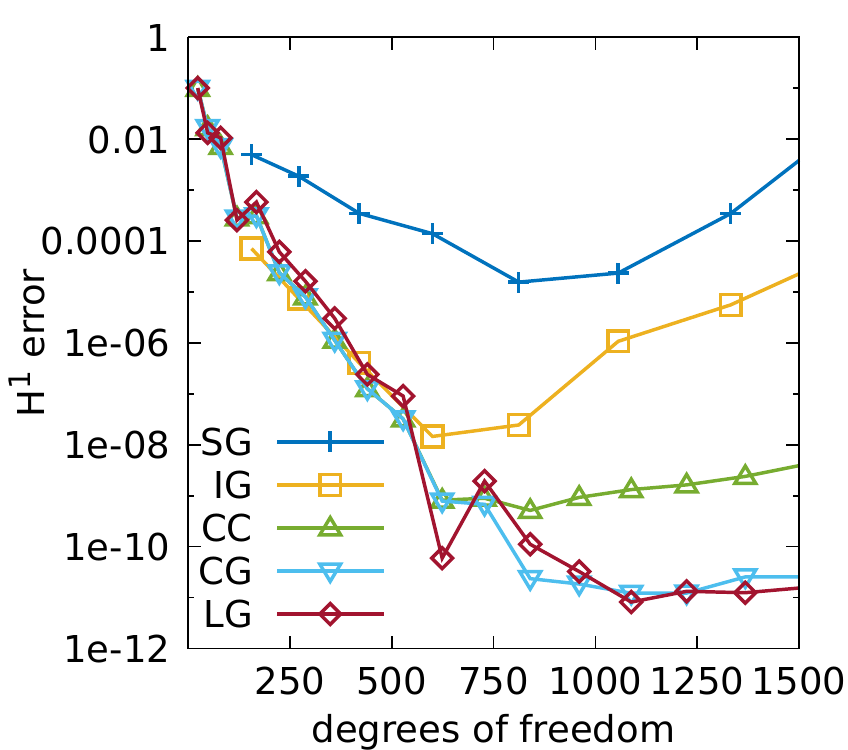}
  }
  \caption{Error after termination of the fixed-point scheme for the Allen--Cahn equation, measured in the $H^1$ norm, as a function of the polynomial degree $p$ (left) and of the estimated operation count vs number of degrees of freedom (right). The figures refer to the annular geometry with forcing term chosen so that the exact solution is $u(\xx)=x_1^2-x_2^3$.}
  \label{fig:allen_cahn}
\end{figure}

The nonlinear test consists in solving the Allen--Cahn equation on the surface shaped as a quarter of annulus, with a forcing term $f$ chosen so that the exact solution is $u(\xx)=x_1^2-x_2^3$.
The tolerance for the fixed point methods, computed as defined in Equation~\eqref{eq:def_increment}, is set to $10^{-15}$.
The results of this test are shown in Figure~\ref{fig:allen_cahn}, where $p$-refinement for the SEMs is compared with $k$-refinement for B-spline and Isogeometric methods. The results confirm what seen in the previous paragraphs: even in nonlinear problems, IGA is extremely efficient for polynomial degrees up to 6 or 8, but fails to achieve spectral accuracy in the context of strong $p$ or $k$-refinement. SEMs, although not competitive with IGA for lower degree polynomials, become the only viable option for polynomial degrees higher than 10. 

\begin{figure}[tbp]
  \centering
  \subfloat[]{
    \includegraphics[width=0.45\textwidth]{Kbeam_log_C1_cond.pdf}
  }
  \subfloat[]{
    \includegraphics[width=0.45\textwidth]{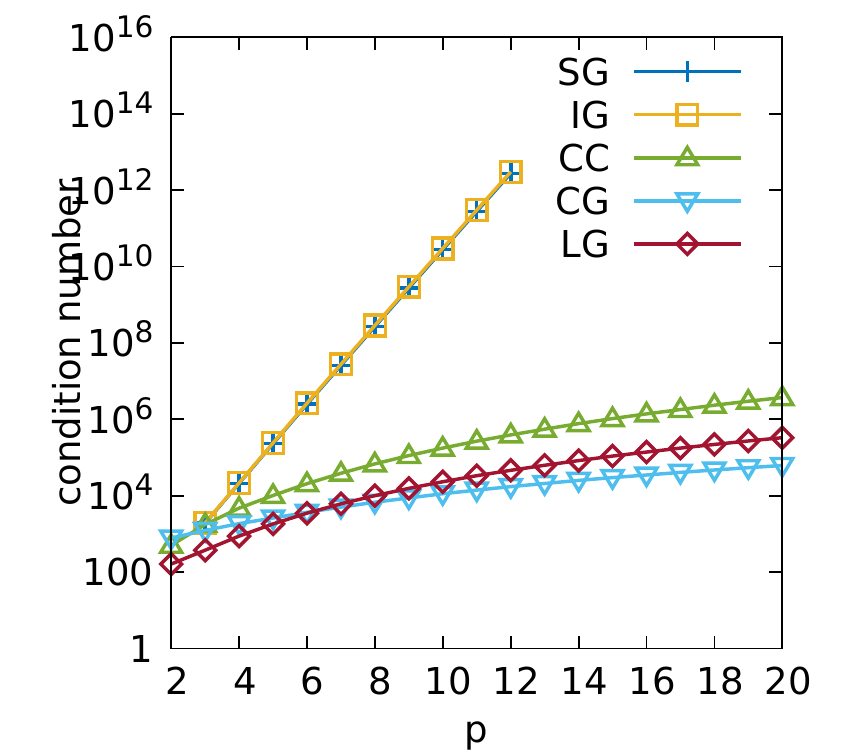}
  } \\
  \subfloat[]{
    \includegraphics[width=0.45\textwidth]{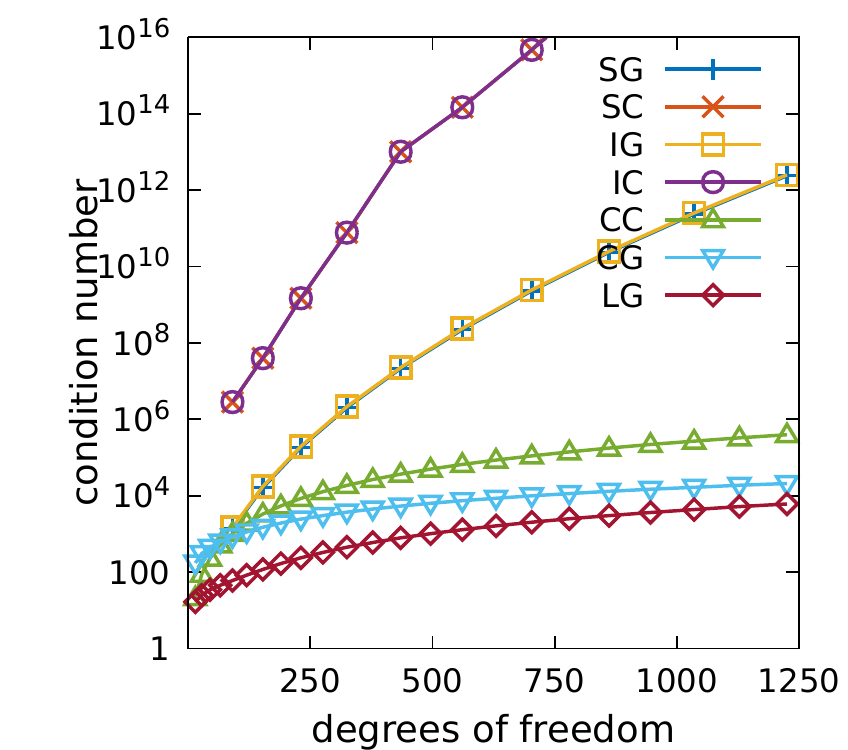}
  }
  \subfloat[]{
    \includegraphics[width=0.45\textwidth]{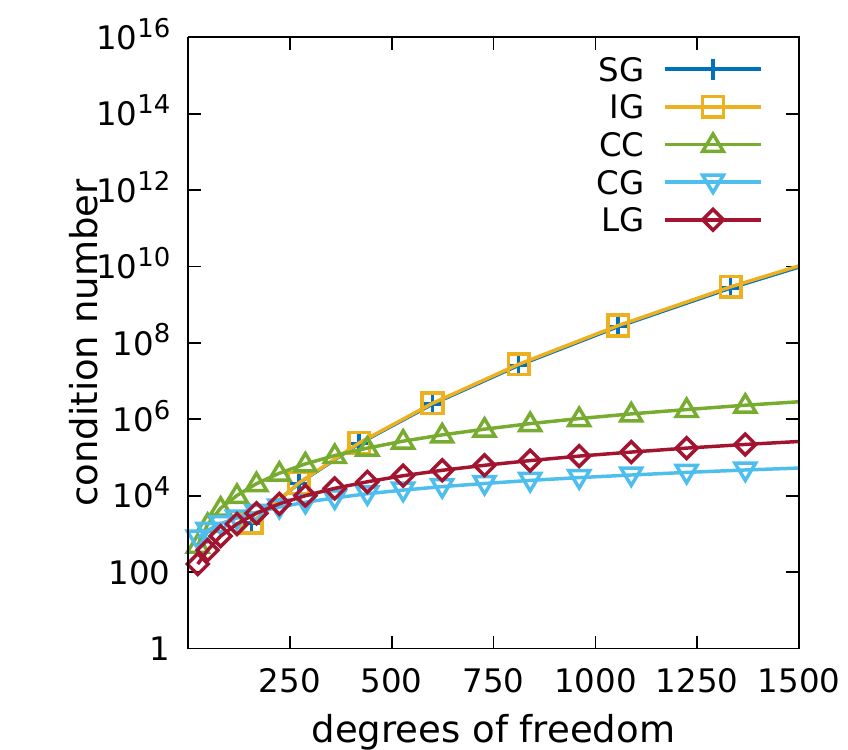}
  }
  \caption{Matrix condition numbers for the $k$-refinement test on the annular surface. On the left, the results refer to the Laplace--Beltrami equation, while on the right figure the results refer to the fixed-point iteration for the Allen--Cahn equation. In the top row, the condition number is plot as a function of the polynomial degree $p$, while in the bottom row the condition number is shown as a function of the number of degrees of freedom.}
  \label{fig:k_refinement_cond}
\end{figure}

For reference, we report the condition number of all matrices, both in the linear and in the nonlinear case, as a function of the polynomial degree and as a function of the number of degrees of freedom in Figure~\ref{fig:k_refinement_cond}. While this comparison is not entirely fair (since the size of the IGA matrices is larger than their spectral counterpart), it is still significant to show the very large rate of growth of the condition number for IGA methods with the polynomial degree. 

\section{Conclusions}
\label{sec:conclusions}
We presented some high-order numerical methods based on NURBS mappings, and applied them to the Laplace--Beltrami equation on some moderately complex surface geometries.

Founding a numerical method on NURBS maps allows to exactly represent the domain geometry and avoids the meshing construction process, a complex and time-consuming step usually done through isoparametric or transfinite maps in $hp$-Finite Element or Spectral Element Methods.

All the numerical methods considered achieve very fast convergence for polynomial degrees between 3 and 8, but only the methods based on Lagrange interpolants at Gauss--Lobatto points reach consistently spectral precision.

The numerical evidence produced in this article raises some limitations to Isogeometric methods in the context of strong $p$-refinement.
In particular, Isogeometric methods perform very well up to a polynomial degree between 8 and 10, then the stiffness matrix becomes too ill-conditioned and starts affecting accuracy.
One notable exception is the Isogeometric Collocation method with globally $\CC^1$ basis functions. This method has reached almost spectral accuracy in all of our numerical tests, at a computational cost comparable with that of the best performer.

Chebyshev multipatch collocation methods composed with NURBS maps seem to be a very good option to achieve high-order and high-precision approximations, avoiding at the same time the need for complex meshing subroutines, at least on two-dimensional problems. The extension to three-dimensional domains would require addressing the many inter-element patching conditions, that would add complexity to the implementation of this method.
An important limitation of Chebyshev and Legendre methods is found in singular meshes. If a side of a quad collapses to a point, the Chebyshev and Legendre methods break due to the loss of degrees of freedom, and to the consequent singularity of the system's matrix. In such a case, B-spline and NURBS methods are more robust and can still provide a solution.

A limitation of Chebyshev methods is that integration by parts can not be carried out due to the presence of a weighted inner product. This can be an issue when higher order differential operators are considered.

\section*{References}
\bibliographystyle{elsarticle-harv}
\bibliography{bibliography}

\end{document}